\documentclass[a4,11pt,epsf,amssymb]{article}
\usepackage{amsmath,amsfonts,latexsym,amssymb, color}
\begin{document}

\title{ \bf Uniform Bahadur Representation for Local Polynomial Estimates of
M-Regression and Its Application to The Additive Model }
\author{ Efang Kong\thanks{%
Eurandom, Technische Universiteit Eindhoven, The Netherlands. %
  E-mail address: \texttt{kong@eurandom.tue.nl.}}\\
%EndAName
{\normalsize \textit{EURANDOM, The Netherlands}}\\
Oliver Linton\thanks{%
Department of Economics, London School of Economics, Houghton
Street, London
WC2A 2AE, United Kingdom. http://econ.lse.ac.uk/staff/olinton/\symbol{126}%
index\_own.html. E-mail address: \texttt{o.linton@lse.ac.uk.}} \\
{\normalsize \textit{London School of Economics, UK }} \\
Yingcun Xia\thanks{%
Department of Statistics and Applied Probability, National University of Singapore, Singapore.
http://www.stat.nus.edu.sg/\symbol{126}%
staxyc.
 E-mail address: \texttt{staxyc@nus.edu.sg.}}\\
{\normalsize \textit{National University of Singapore, Singapore}} }
\date{}
\def\beginn{\begin{eqnarray*}}
\def\endn{\end{eqnarray*}}
\def\beginy{\begin{eqnarray}}
\def\endy{\end{eqnarray}}
\def\n{\nonumber}
\newtheorem{Theorem}{Theorem}[section]
\newtheorem{Example}[Theorem]{Example}
\newtheorem{Lemma}[Theorem]{Lemma}
\newtheorem{Note}[Theorem]{Note}
\newtheorem{Proposition}[Theorem]{Proposition}
\newtheorem{Corollary}[Theorem]{Corollary}
\newtheorem{Remark}[Theorem]{Remark}
\newcommand{\un}[1]{\underline{#1}}

\def\X{{\mathbf X}}
\def\x{\textit{x}}
\def\Cov{\mbox{\rm Cov}}
\def\Var{\mbox{\rm Var}}
\def\YX{_{Y|X}}
\def\argmin{\mbox{\rm arg}\min}
\def\b{{\mathbf b}}
\def\e{{\varepsilon}}
\def\F{{\mathbf F}}
\def\be{\mbox{ e}}
\def\bu{\mbox{ U}}
\def\bv{\mbox{ V}}
\def\CV{\mbox{\small CV}}
\def\cov{\mbox{ cov}}
\def\T{\mbox{\rm T}}
\def\L{\mbox{\rm L}}
\def\argmax{\mbox{\rm argmax}}

\def\pperp{\perp\hspace{-.25cm}\perp}
\def\t{{\hspace{-0.05cm}\top}}
\def\A{{\cal A}}
\def\Bt{\theta^\top}
\def\E{{\cal E}}
\def\D{{\cal D}}
\def\N{{\cal N}}
\def\M{{{\cal M}}}
\def\R{{\cal R}}
\def\rt{\raisebox{1.5ex}[0pt]}
\def\btd{\bigtriangledown}
\def\btdt{\bigtriangledown^\t\hspace{-0.1cm}}
\def\dfor{{\qquad \mbox{\rm for} \quad}}
\def\C{{\cal C}}
\def\g{{\sl g}}
\def\rt{\raisebox{1.8ex}[0pt]}
\def\rtf{\raisebox{1.2ex}[0pt]}
\def\ditem{\vspace{-0.2cm} \item}

\maketitle

\begin{center}{\large S}UMMARY\end{center}
\baselineskip1.9em {\small We use local polynomial fitting to
estimate the nonparametric M-regression function for strongly
mixing stationary processes $\{(Y_{i},\underline{X}_{i})\}$. We
establish a strong uniform consistency rate for the Bahadur
representation of estimators of the regression function and its
derivatives. These results are fundamental for statistical
inference and for applications that involve plugging in such
estimators into other functionals where some control over higher
order terms are required. We apply our results to the estimation
of an additive M-regression model. }


\baselineskip2.1em

\def\x{{\bf x}}

\section{Introduction}

In many contexts one wants to evaluate the properties of some
procedure that is a functional of some given estimators. It is
useful to be able to work with some plausible high level assumptions
about those estimators rather than to rederive their properties for
each different application. In a fully parametric (and stationary,
weakly dependent data) context it is quite common to assume that
estimators are root-n consistent and asymptotically normal. In some
cases this property suffices; in other cases one needs to be more
explicit in terms of the linear expansion of these estimators, but
in any case such expansions are quite natural and widely applicable.
In a nonparametric context there is less agreement about the use of
such expansions and one often sees standard properties of standard
estimators derived anew for a different purpose. It is our objective
to provide results that can circumvent this. The types of
application we have in mind are estimation of semiparametric models
where the parameters of interest are explicit or implicit
functionals of nonparametric regression functions and their
derivatives; see Powell (1994), Andrews (1994) and Chen, Linton and
Van Keilegom (2003). Another class of applications includes
estimation of structured nonparametric models like additive models
(Linton and Nielsen, 1995) or generalized additive models (Linton,
Sperlich and Van Keilegom, 2007).

We motivate our results in a simple i.i.d. setting. Suppose we have
a random sample $\{Y_{i},X_{i}\}_{i=1}^{n}$ and consider the
Nadaraya-Watson estimator
of the regression function $m(x)=E(Y_{i}|X_{i}=x),$%
\beginn
\hat{m}(x)=\frac{\hat{r}(x)}{\hat{f}(x)}=\frac{n^{-1}\sum
_{i=1}^{n}K_{h}(X_{i}-x)Y_{i}}{n^{-1}\sum_{i=1}^{n}K_{h}(X_{i}-x)},
\endn
where $K$ is a symmetric density function, $h$ is a bandwidth and
$K_{h}(.)=K(./h)/h.$ Standard arguments (H\"{a}rdle, 1990) show that
under suitable smoothness conditions
\beginy
\hat{m}(x)-m(x)=h^{2}b(x)+\frac{1}{nf(x)}\sum_{i=1}^{n}K_{h}%
(X_{i}-x)\varepsilon_{i}+R_{n}(x),\label{bo}
\endy
where $b(x)=\int u^{2}K(u)du[m^{\prime \prime }(x)+2m^{\prime
}(x)f^{\prime }(x)/f(x)]/2,$ while $f(x)$ is the covariate density
and $\varepsilon _{i}\equiv Y_{i}-m(X_{i}) $ is the error term. The
remainder term $R_{n}(x)$ is of smaller order (almost surely) than
the two leading terms. Such an expansion is sufficient to derive the
central limit theorem for $\hat{m}(x)$ itself, but generally is not
 sufficient if $\hat{m}(x)$ is to be plugged into some semiparametric
procedure. For example, suppose we estimate the
parameter $\theta _{0}=\int m(x)^{2}dx$ by $\hat{\theta}=\int \hat{m}%
(x)^{2}dx$, where the integral is over some compact set
${\mathcal{D}};$ we would expect to find that
$n^{1/2}(\hat{\theta}-\theta _{0})$ is asymptotically normal. Based
on expansion (\ref{bo}), the argument goes like this. First, we
obtain the following
\beginn
n^{1/2}(\hat{\theta}-\theta_{0})=2n^{1/2}\int m(x)\{\hat
{m}(x)-m(x)\}dx+n^{1/2}\int[\hat{m}(x)-m(x)]^{2}dx.
\endn
If it can be shown that $\hat{m}(x)-m(x)=o(n^{-1/4})$ a.s. uniformly
in $x\in {\mathcal D}$ ( such results are widely available; see for
example Masry (1996)), we have
\beginn
n^{1/2}(\hat{\theta }-\theta _{0})=2n^{1/2}\int m(x)\{\hat{m}%
(x)-m(x)\}dx+o(1),\quad a.s.
\endn
Note that the quantity on the right hand side is the term in
assumption 2.6 of Chen, Linton, and Van Keilegom (2003) which is
assumed to be asymptotically normal. It is the verification of this
condition with which we are now
concerned. If we substitute in the expansion (\ref{bo}) we obtain%
\begin{align*}
n^{1/2}(\hat{\theta }-\theta _{0})& =2n^{1/2}h^{2}\int
m(x)b(x)dx+2n^{1/2}\int \frac{m(x)}{f(x)}n^{-1}\sum_{i=1}^{n}K_{h}(X_{i}-x)%
\varepsilon _{i}dx \\
& \qquad +2n^{1/2}\int m(x)R_{n}(x)dx+o(1),\quad a.s.
\end{align*}
If $nh^{4}\rightarrow0,$ then the first term (the smoothing bias
term) is $o(1).$ By a change of variable, the second term (the
stochastic term) can be written as a sum of  independent random
variables with zero mean
\beginn
&n^{1/2}\int
m(x)f^{-1}(x)n^{-1}\sum\limits_{i=1}^{n}K_{h}(X_{i}-x)\varepsilon
_{i}dx=n^{-1/2}\sum\limits_{i=1}^{n}\xi_{n}(X_{i})\varepsilon_{i},&\hspace{.2cm}\\
&\xi_{n}(X_{i})=\int m(X_{i}+uh)f^{-1}(X_{i}+uh)K(u)du,&
\endn
and this term obeys the Lindeberg central limit theorem under
standard conditions. The problem is that equation (\ref{bo}) only
guarantees that $\int m(x)R_{n}(x)dx=o(n^{-2/5})$ a.s. at best$.$
Actually, in this case it is possible to derive a more useful
Bahadur expansion (Bahadur, 1966) for the kernel estimator
\begin{equation}
\hat{m}(x)-m(x)=h^{2}b_{n}(x)+\{E\hat{f}(x)\}^{-1}n^{-1}\sum
_{i=1}^{n}K_{h}(X_{i}-x)\varepsilon_{i}+R_{n}^{\ast}(x),\label{boo}%
\end{equation}
where $b_{n}(x)$ is deterministic and satisfies $b_{n}(x)\rightarrow
b(x)$ uniformly in $x\in {\mathcal D},$ and $E\hat{f}(x)\rightarrow
f(x)$ uniformly in
$x\in {\mathcal D},$ while the remainder term now satisfies%
\begin{equation}
\sup_{x\in {\mathcal D}}\left\vert R_{n}^{\ast}(x)\right\vert =O\left( \frac{\log n}{nh%
}\right) \quad a.s.\label{b1}%
\end{equation}
This property is a consequence of the uniform convergence rate of $%
\hat{f}(x)-E\hat{f}(x),$ $n^{-1}\sum_{i=1}^{n}K_{h}(x$ $-X_{i})%
\{m(X_{i})-m(x)\}-EK_{h}(X_{i}-x)\{m(X_{i})-m(x)\}$ and $%
n^{-1}\sum_{i=1}^{n}K_{h}(X_{i}-x)\varepsilon _{i}$ that follow
from, for example Masry (1996). Clearly, by appropriate
choice of $h$, $R_{n}^{\ast }(x)$ can be made to be $%
o(n^{-1/2})$ $a.s.$ uniformly over ${\mathcal D}$ and thus
$2n^{1/2}\int m(x)R_{n}^{\ast }(x)dx=o(1)$ $a.s.$. Therefore, to
derive asymptotic normality for $n^{1/2}(\hat{\theta }-\theta
_{0})$, one can just work with the two leading terms in (\ref{boo}).
These terms are slightly more complicated than in the previous
expansion but are still sufficiently simple for many purposes; in
particular, $b_{n}(x)$ is uniformly bounded so that
provided $nh^{4}\rightarrow 0,$ the smoothing bias term satisfies $h^{2}%
n^{1/2}\int m(x)b_{n}(x)dx\rightarrow 0,$ while the stochastic term
is a sum of zero mean independent random variables
\begin{equation*}
n^{1/2}\int \frac{m(x)}{\overline{f}(x)}n^{-1}\sum_{i=1}^{n}K_{h}(X_{i}-x)%
\varepsilon _{i}dx=n^{-1/2}\sum_{i=1}^{n}\overline{\xi }%
_{n}(X_{i})\varepsilon _{i}
\end{equation*}%
\begin{equation*}
\overline{\xi }_{n}(X_{i})=\int \frac{m(X_{i}+uh)}{\overline{f}(X_{i}+uh)}%
K(u)du,
\end{equation*}%
and obeys the Lindeberg central limit theorem under standard
conditions, where $\overline{f}(x)=E\hat{f}(x).$ This argument shows
the utility of the Bahadur expansion (\ref{boo}). There are many
other applications of this result because a host of probabilistic
results are available for random variables like
$n^{-1}\sum_{i=1}^{n}K_{h}(X_{i}-x)\varepsilon _{i}$ and integrals
thereof$.$

The one-dimensional Nadaraya-Watson estimator for i.i.d. data is
particularly easy to analyze and the above arguments are well known.
However, the limitations of this estimator are manyfold and there
are good theoretical reasons for working instead with the local
polynomial class of estimators (Fan and Gijbels, 1996). In addition,
for many data especially financial time series data one may have
concerns about heavy tails or outliers that point in the direction
of using robust estimators like the local median or local quantile
method, perhaps combined with local polynomial fitting. We examine a
general class of (nonlinear) M-regression function (that is,
location functionals defined through minimization of a general
objective function $\rho(.)$) and derivative estimators. We treat a
general time series setting where the multivariate data are strongly
mixing. Under mild conditions, we establish a uniform strong Bahadur
expansion like (\ref{boo}) and (\ref{b1}) with remainder term of
order $(\log n/nh^{d})^{3/4}$ almost surely, which is almost optimal
or in other words can't be improved further based on the results in
Kiefer (1967) under i.i.d. setting. The leading terms are linear and
functionals of them can be analyzed simply. The remainder term can
be made to be $o(n^{-1/2})$ a.s. under restrictions on the
dimensionality in relation to the amount of smoothness possessed by
the M-regression function.

The best convergence rate of unrestricted nonparametric estimators
strongly depends on $d$, the dimension of $\underline{x}$. The rate
decreases dramatically as $d$ increases (Stone, 1982). This
phenomenon is the so-called \textquotedblleft curse of
dimensionality\textquotedblright . One approach to reduce the curse
is by imposing model structure. A popular model structure is the
additive model assuming that
\beginy
m(x_{1},\ldots ,x_{d})=c+m_{1}(x_{1})+...+m_{d}(x_{d}),
\label{junk}
\endy
where $c$ is an unknown constant and $m_{k}(.),\ k=1,\ldots ,d$ are
unknown functions which have been normalized such that
$Em_{k}(\mathbf{x}_{k})=0$ for $k=1,\ldots ,d.$ In this case, the
optimal rate of convergence is the same as in univariate
nonparametric regression (Stone, 1986). An additive M-regression
function is given by (\ref{junk}) with $m(x)$ being the M-regression
function defined in (\ref{def0}). Previous work on additive quantile
regression, for example, includes Linton (2001) and Horowitz and Lee
(2005) for the i.i.d. case. An interesting application of the
additive M-regression model is to combine (\ref{junk}) with the
volatility model
\begin{equation*}
Y_{i}=\sigma _{i}\varepsilon _{i}\quad \mbox{and}\quad \ln \sigma
_{i}^{2}=m(X_{i}),
\end{equation*}%
where $X_{i}=(Y_{i-1},\ldots ,Y_{i-d})^{\top }.$ We suppose that $%
\varepsilon _{i}$ satisfies $E[\varphi (\ln \varepsilon
_{i}^{2};0)|X_{i}]=0
$ for some function $\varphi (.)$, whence $m(.)$ is the conditional $M$-regression of $\ln Y_{i}^{2}$ given $%
X_{i}.$ Peng and Yao (2003) have applied LAD estimation to
parametric ARCH and GARCH models and have shown the superior
robustness property of this procedure over Gaussian QMLE with regard
to heavy tailed innovations. This heavy tail issue also arises in
nonparametric regression models, which is why our procedures may be
useful. Empirical evidence also suggest that moderately high
frequency financial data are often heavy tailed. We apply our
Bahadur expansions to the study of marginal integration estimators
(Linton and Nielsen, 1995) of the component functions in additive
M-regression model in which
case we only need the remainder term to be $o(n^{-p/(2p+1)})$ a.s.$,$ where $%
p$ is a smoothness index.

Bahadur representations (Bahadur, 1966) have been widely studied and
applied, with notable refinements in the i.i.d. setting by Kiefer
(1967). A recent paper of Wu (2005) extends these results to a
general class of dependent processes and provides a review. The
closest paper to ours is Hong (2003) who established a Bahadur
representation for essentially the same local polynomial
M-regression estimator as ours. However, his results are: (a)
pointwise, i.e., for a single $x$ only; (b) the covariates are
univariate; (c) for i.i.d. data. Clearly, this limits the range of
applicability of his results, and specifically, the applications to
semiparametric or additive models are perforce precluded.

\section{The General Setting}

Let $\{(Y_{i},\un{X}_i)\}$ be a jointly stationary processes, where
$\un{X}_i= (\x_{i1}, ..., \x_{id})^\t$ with $d\ge 1$ and $Y_{i}$ is
a scalar. As dependent observations are considered
in this paper, we introduce here the mixing coefficient. Let ${\mathbf{F}}%
_{s}^{t}$ be the $\sigma -$ algebra of events generated by random variables $%
\{(Y_{i},\underline{X}_{i}),s\leq i\leq t\}.$ A stationary stochastic processes $%
\{(Y_{i},\underline{X}_{i})\}$ is strongly mixing if
\begin{equation*}
\sup\limits_{\overset{A\in {\mathbf{F}}_{-\infty }^{0}}{B\in {\mathbf{F}}%
_{k}^{\infty }}}|P[AB]-P[A]P[B]|=\gamma \lbrack k]\rightarrow 0,\ \mbox{as }%
k\rightarrow \infty ,
\end{equation*}%
and $\gamma \lbrack k]$ is called the strong mixing
coefficient.

Suppose $\rho(.;.)$ is a loss function. Our first goal is to
estimate
 the multivariate  M-regression function
 \beginy
 m(x_1,\cdots,x_d)=\argmin\limits_{\theta}
 E\{\rho(Y_i;\theta)|\un{X}_i=(x_1,\cdots,x_d)\},\label{def0}
 \endy
and its partial derivatives based on  observations
$\{(Y_{i},\un{X}_i)\}_{i=1}^n$. An important  example of the
M-function is the $q-$th ($0<q<1$) quantile of $Y_i$ given
$\un{X}_i=(x_1,\cdots,x_d)^\t$, with loss function given by
$\rho(y;\theta)=(2q-1)(y-\theta)+|y-\theta|$. Another example is the
$L_{q}$ criterion $\rho (y;\theta )=|y-\theta |^{q}$ for $q>1,$
which includes the least square criterion $\rho (y;\theta
)=(y-\theta )^{2}$ with  $m(.)$  the conditional expectation
 of $Y_{i}$ given $\underline{X}_{i}.$ Yet another  example is the celebrated
Huber's function (Huber, 1973)
\beginy
\rho(t)=t^2/2I\{|t|<k\}+(k|t|-k^2/2)I\{|t|\ge k\}.\label{Huber}
\endy
Suppose $m(\un{x})$ is differentiable up to order $p+1$ at $\un{x}=
(x_1, ..., x_d)^\t$. Then  the multivariate $p$'th order local
polynomial approximation of $m(\un{z})$ for any $ \un{z}$ close to $
\un{x}$ is given by
\beginn
m(\un{z})=\sum\limits_{0\le |{\un{r}}|\le
p}\frac{1}{\un{r}!}D^{\un{r}}m(\un{x})(\un{z}-\un{x})^{\un{r}},
\endn
where $ {\un{r}} = (r_1, ..., r_d),\  |{\un{r}}| = \sum_{i=1}^d
r_i,\ \un{r}!= r_1!\times\cdots\times r_d!$ and
\beginy
D^{\un{r}}m(\un{x})=\frac{\partial^{\un{r}}m(\un{x})}{\partial
x_1^{r_1}\cdots\partial x_d^{r_d}},\quad \un{x}^{\un{r}} =
x_1^{r_1} \times
 ... \times x_d^{r_d},
 \quad \sum_{ 0 \le |{\un{r}}| \le p} =
 \sum_{j=0}^p
\underset{r_1+...+r_d=j}{\sum_{r_1=0}^j ... \sum_{r_d=0}^j}.
\endy
Let $K(\un{u})$ be a density function on $R^d$, $h$  a bandwidth and
$K_h(u)=K(u/h)$. With  observations $\{(Y_{i},\un{X}_i)\}_{i=1}^n$,
we consider minimizing the following quantity with respect to
$\beta_{\un{r}},\ 0 \le |\un{r}| \le p $
\beginy
\sum\limits_{i=1}^{n}K_h(\un{X}_i-\un{x})\rho\Big(Y_{i};\sum_{0 \le |{\un{r}}| \le
p} \beta_{\un{r}}(\un{X}_i-\un{x})^{\un{r}}\Big).\label{mini}
\endy
Denote by $\hat \beta_{\un{r}}(\un{x}),\ 0 \le |r| \le p, $ the
minima of (\ref{mini}). The M-function $m(\un{x}) $ and its
derivatives $D^{\un{r}}m(\un{x})$ are then estimated respectively by
\beginy
\hat m(\un{x}) = \hat \beta_{\un{0}}(\un{x}) \quad  \mbox{and} \quad
\hat D^{\un{r}}m(\un{x}) = \un{r}!\hat \beta_{\un{r}}(\un{x}),\ 1
\le |{\un{r}}| \le p.\label{why}
\endy

\section{Main Results}
In Theorem 3.2 below we give our main result, the uniform strong
Bahadur representation for the vector
$\hat{\beta}_{p}(\underline{x}).$ We first need to develop some
notations to define the leading terms in the expansion.

Let $N_i={i+d-1\choose d-1}$ be the number of distinct $d-$tuples
$\un{r}$ with $|{\un{r}}|=i$. Arrange these $ d-$tuples as a
sequence in a lexicographical order(with the highest priority given
to the last position so that $(0,\cdots,0,i)$ is the first element
in the sequence and $(i,0,\cdots,0)$ the last element).  Let
$\tau_i$
 denote this $1$-to-$1$ mapping, i.e. $\tau_i(1)=(0,\cdots,0,i),\cdots,
 \tau_i(N_i)=(i,0,\cdots,0).$
 For each $i=1,\cdots,p$, define a  $N_{i}\times 1$
vector $\mu_{i}(\un{x})$ with its $k$th element given by
$\un{x}^{\tau_{i}(k)}$ and write
$\mu(\un{x})=(1,\mu_{1}(\un{x})^\t,\cdots,\mu_{p}(\un{x})^\t)^\t, $
which is a column vector of length $N=\sum_{i=0}^pN_i.$ Similarly
define  vectors $\beta_p(\un{x})$ and $\un{\beta}$ through the same
lexicographical arrangement of $D^{\un{r}}m(\un{x})$ and
$\beta_{\un{r}}$ in (\ref{mini}) for $0 \le |{\un{r}}| \le p.$ Thus
(\ref{mini}) can be rewritten as
\beginy
\sum\limits_{i=1}^{n}K_h(\un{X}_i-\un{x})\rho(Y_{i};\mu(\un{X}_i-\un{x})^\t\un{\beta}).\label{mini2}
\endy
Suppose the minimizer of (\ref{mini2}) is denoted as $\tilde
\beta_n(\un{x})$. Let $\hat \beta_p(\un{x})=W_p\hat\beta_n(\un{x}),$
 where $W_p$ is the diagonal matrix with diagonal entries
the lexicographical arrangement of $\un{r}!,\ 0 \le |{\un{r}}| \le
p.$

Let $\nu_{\un{i}}=\int K(\un{u})\un{u}^{\un{i}}d\un{u}$. For $g(.)$
given in  (\ref{aa3}), define
\beginn
\nu_{n\un{i}}(\un{x})=\int
K(\un{u})\un{u}^{\un{i}}g(\un{x}+h\un{u})
f(\un{x}+h\un{u})d\un{u}.
\endn
For $0\le j,k\le p$, let $S_{j,k}$ and $S_{n,j,k}(\un{x})$ be two
$N_j\times N_k$ matrices   with  their $(l,m)$ elements respectively
given by
\beginy
\Big[S_{j,k}\Big]_{l,m}=\nu_{\tau_{j}(l)+\tau_{k}(m)}(\un{x}),\quad
\Big[S_{n,j,k}(\un{x})\Big]_{l,m}=\nu_{n,\tau_{j}(l)+\tau_{k}(m)}(\un{x}).\label{von}
\endy
Now define the $N\times N$ matrices $S_p$ and $S_{n,p}(\un{x})$ by
\beginn
S_p=\left[
                 \begin{array}{cccc}
                  S_{ 0,0}&S_{ 0,1}&\cdots& S_{ 0,p}\\
                 S_{ 1,0}&S_{ 1,1}&\cdots& S_{ 1,p}\\
                 \vdots&\ddots&\vdots\\
                   S_{ p,0}&S_{ p,1}&\cdots& S_{ p,p}\\
                 \end{array}
               \right],
\quad
S_{n,p}(\un{x})=\left[
                 \begin{array}{cccc}
                  S_{n,0,0}(\un{x})&S_{n,0,1}(\un{x})&\cdots& S_{n,0,p}(\un{x})\\
                 S_{n,1,0}(\un{x})&S_{n,1,1}(\un{x})&\cdots& S_{n,1,p}(\un{x})\\
                 \vdots&\ddots&\vdots\\
                   S_{n,p,0}(\un{x})&S_{n,p,1}(\un{x})&\cdots& S_{n,p,p}(\un{x})\\
                 \end{array}
               \right].
\endn
According to Lemma \ref{Masry},  $S_{n,p}(\un{x})$ converges to
$g(\un{x})f(\un{x})S_p$ uniformly in $\un{x}\in {\mathcal D}$ almost
surely. Hence for $|S_p|\ne 0$, we can define
\beginy
\beta^*_{n}(\un{x})=-\frac{1}{nh^d}W_pS_{n,p}^{-1}(\un{x})H_n^{-1}\sum\limits_{i=1}^{n}
K_h(\un{X}_i-\un{x})\varphi(Y_{i},\mu(\un{X}_i-\un{x})^\t
\beta_p(\un{x}))\mu(\un{X}_i-\un{x}), \label{matrix}
\endy
where $\varphi(.;.)$ is the piecewise derivative of $\rho(.,.)$, as
defined in (A1) and $H_n$ is the diagonal matrix with diagonal
entries $h^ {|\un{r}|},\ 0 \le |{\un{r}}| \le p$ in the
aforementioned lexicographical order. The quantity $\beta _{n}^{\ast
}(\underline{x})$ is the leading term of our expansion; it contains
both a bias term, $E\beta _{n}^{\ast }(\underline{x}),$ and a
stochastic leading term $\beta _{n}^{\ast }(\underline{x})-E\beta
_{n}^{\ast }(\underline{x})$.

Denote the typical element of $\beta^*_{n}(\un{x})$   by
$\beta^*_{n\un{r}}(\un{x}),\ 0 \le |{\un{r}}| \le p$ and the density
function of $\un{X}$ by $f(.)$. The following results on
$E\beta^*_{n\un{r}}(\un{x})$ is an extension of Proposition 2.2 in
Hong (2003) to the multivariate case.
\begin{Proposition}\label{Pro1}
If $f(\un{x})>0$
  and conditions (A1)-(A5) in the Appendix hold, then
\beginn
E\beta^*_{n\un{r}}(\un{x})=\left\{
\begin{array}{ll}
-h^{p+1}e_{N(\un{r})}W_pS_{p}^{-1}B_{1}{\bf
m}_{p+1}(\un{x})+o(h^{p+1}), &
 \hspace{-2.5cm}\mbox{ for } p-|\un{r}|\mbox{ odd,}\\
 \\
 -h^{p+2}e_{N(\un{r})}W_pS_{p}^{-1}\Big[\{f g\}^{-1}(\un{x}){\bf
m}_{p+1}(\un{x})\{\tilde M(\un{x})-N_pS_{p}^{-1}B_{1}\}+B_{2}{\bf
m}_{p+2}(\un{x})\Big]\\
+o(h^{p+2}), & \hspace{-2.5cm} \mbox{ for }p-|\un{r}|\mbox{ even,}
\end{array}
\right.
\endn
where
$N(\un{r})=\tau_{|\un{r}|}^{-1}(\un{r})+\sum_{k=0}^{|\un{r}|-1}N_k$,
$e_i$ is a $N\times 1$ vector having $1$ as the $i$th entry with all
other entries $0,$ and $ B_{1}=\left[S_{0,p+1}, S_{1,p+1}, \cdots
S_{p,p+1} \right]^\t,\ B_{2}=\left[S_{0,p+2}, S_{1,p+2}, \cdots
S_{p,p+2} \right]^\t. $
\end{Proposition}
We next present our main result, the Bahadur representation for
local polynomial estimates $\hat \beta_p(\un{x})$.
\begin{Theorem}
\label{T1} Suppose  (A1)-(A7) in the Appendix hold with
$\lambda_2=(p+1)/2(p+s+1)$ for some $s\ge 0$ and  ${\mathcal D}$ is
any compact subset of $R^d.$ Then
\beginn
\sup\limits_{\un{x}\in {\mathcal D} }|H_n\{\hat
\beta_p(\un{x})-\beta_p(\un{x})\}-\beta_{n}^*(\un{x})|=
O\Big(\Big\{\frac{\log n}{nh^d}\Big\}^{\lambda(s)}\Big) \
\mbox{almost surely},\label{mathisca}
\endn
where  $|.|$ is taken to be the sup norm  and
\beginn
\lambda(s)=\min\Big\{\frac{p+1}{p+s+1},\ \frac{3p+3+2s}{4p+4s+4}\Big\}.
\endn
\end{Theorem}
\noindent{\bf Remark 1}.  According to Theorem 1 in Kiefer (1967),
 the point-wise sharpest bound of the remainder term in Bahadur representation of the sample
quantiles is $(\log\log n/n)^{3/4}$.
 As $\lambda(0)=3/4$,  we could safely claim  the results here could not be
further improved for a general class of loss functions $\rho(.)$
specified by (A1) and (A2).  Nevertheless, it is possible to derive
stronger results, if the concerned loss functions enjoy higher
degree of smoothness; see (\ref{b1}) in which case $\rho(.)$ is the
squared loss function. More specifically, suppose $\varphi(.)$ is
Lipschitz continuous and (A1)-(A7) in the Appendix hold with
$\lambda_2=1/2$ and $\lambda_1=1$.  Then we  prove in the Appendix
that with probability $1$ and uniformly in $\un{x}\in {\mathcal D}$,
\beginy
\sup\limits_{\un{x}\in {\mathcal D} }|H_n\{\hat
\beta_p(\un{x})-\beta_p(\un{x})\} -\beta_{n}^*(\un{x})|
=O\Big(\frac{\log n}{nh^d}\Big) \mbox{ almost
surely}.\label{mathisca}
\endy

\noindent{\bf Remark 2}. The dependence among the observations
doesn't have any impact on  the rate of
 uniform convergence, given that the degree of the dependence, as measured by
 the mixing coefficient
 $\gamma[k]$, is weak enough such that (\ref{centraal}) and (\ref{rbig2}) are satisfied. This is in accordance with
 the results in Masry (1996), where he proved that for local polynomial estimator of the
 conditional mean function, the uniform convergence rate  is
 $(nh^d/\log n)^{-1/2}$, the same as in the independent case.

 \noindent{\bf Remark 3}. It is of practical interest
to provide an explicit rate of decay for the strong mixing
coefficient $\gamma[k]$ of the form $\gamma[k]=O(1/k^c)$ for some
$c>0$ (to be determined) for Theorem \ref{T1} to hold. It is easily
seen that, among all the conditions imposed on $\gamma[k]$, the
summability condition (\ref{rbig2}) is the most restrictive. We
assume that
\beginn
h=h_n\sim (\log n/n)^{\bar a} \mbox{ for some }
\frac{1}{2(p+s+1)+d}\le \bar a<
\frac{1}{d}\Big\{1-\frac{4}{(1-\lambda_2)\nu_2
-4\lambda_1+2(1+\lambda_2)}\Big\}
\endn
whence (\ref{rbig1}) is satisfied. Algebraic calculations show that
the summability condition (\ref{rbig2}) is satisfied provided that
\beginy
c>\nu_2\frac{(1-\bar a
d)\{(1-\lambda_2)(4N+1)+8N\lambda_1\}+10+(4+8N)\bar a d}
{2(1-\lambda_2)(1-\bar a d)\nu_2-8\bar a d+4(1-\bar a
d)(1-\lambda_2-2\lambda_1)}-1\equiv c(d,p,\nu_2,\bar
a,\lambda_1,\lambda_2). \label{pot}
\endy
Note that we would need the following condition
\beginn
\nu_2>2+\frac{4\{\bar a d+(1-\bar a d)\lambda_1\}}{(1-\bar a
d)(1-\lambda_2)}
\endn
 to secure  positive denominator for (\ref{pot}). As
  $c(d,p,\nu_2,\bar a,$ $\lambda_1,\lambda_2)$ is decreasing in
$\nu_2(\le \nu_1)$, there is a tradeoff between the order $\nu_1$ of
the moment $E|\varphi(\e_i)|^{\nu_1}<\infty$ and the decay rate of
the strong mixing coefficient $\gamma[k]$: the existence of higher
order moments allows $\gamma[k]$ to decay more slowly.

\noindent{\bf Remark 4}. It is trivial to generalize the result in
Theorem 3.2 to functionals of the M-estimates
$\hat\beta_{p}(\un{x})$. Denote the typical elements of
$\hat\beta_p(\un{x})\ $ and $\beta_p(\un{x})$ by
$\hat\beta_{p\un{r}}(\un{x})$ and $\beta_{p\un{r}}(\un{x}),\ 0\le
|\un{r}|\le p$ respectively. Suppose $G(.): \ R^d\to R$ satisfies
that for any compact set $\mathcal D\subset R^d$, there exists some
constant $C>0$,
 such that
 $|G',(\beta_{p\un{r}}(\un{x}))|\le C$ and $|G{''}(\beta_{p\un{r}}(\un{x}))|\le C $
  for all $\un{x}\in {\mathcal D}$. Then  with probability $1$,
\beginy
\sup\limits_{\un{x}\in {\mathcal D}
}\Big|h^{|\un{r}|}\left[G\{\hat\beta_{p\un{r}}(\un{x})\}-G\{\beta_{p\un{r}}(\un{x})\right]
-G'\{\beta_{p\un{r}}(\un{x})\} \beta^*_{n\un{r}}(\un{x})\Big|=
O\Big(\Big\{\frac{\log n}{nh^d}\Big\}^{\lambda(s)}\Big) \label{gexp}
\endy
uniformly for all $\un{x}\in {\mathcal D}$.

The following proposition follows from Theorem \ref{T1}  and uniform
convergence of sum of weakly dependent zero mean random variables.

\begin{Corollary}\label{Pro2}Suppose
conditions in Theorem \ref{T1} hold with $s=0$.   Then with
probability $1$ we have, uniformly in $\un{x}\in {\mathcal D}$,
\beginn
H_n\{\hat \beta_p(\un{x})-\beta_p(\un{x})\} -E\beta_{n}^*(\un{x})-
\frac{W_pH_n^{-1}}{nh^d}S_{np}^{-1}(\un{x}) \sum\limits_{i=1}^{n}
K_h(\un{X}_i-\un{x})\varphi(\varepsilon_{i})\mu(\un{X}_i-\un{x})
=O\Big(\Big\{\frac{\log n}{nh^d}\Big\}^{3/4}\Big).
\endn

\end{Corollary}

\section{M-Estimation of the Additive model}

In this section, we apply our main result to derive the properties
of a class of estimators in the additive M-regression model
(\ref{junk}). In terms of estimating the component functions
$m_{k}(.),\ k=1,\ldots ,d$ in (\ref{junk}),  the marginal
integration method (Linton and Nielsen, 1995)  is known to achieve
the optimal rate under certain conditions. This involves estimating
first the unrestricted M-regression function $m(.)$ and then
integrating it over some directions.
Partition $\underline{X}_{i}=(x_{1},\ldots ,x_{d})$ as $\underline{X}_{i}=(%
\mathbf{x}_{1i},\underline{X}_{2i})$, where $\mathbf{x}_{1i}$ is the
one dimensional direction of interest and $\underline{X}_{2i}$ is a
$d-1$ dimensional nuisance direction. Let
$\underline{x}=(x_{1},\underline{x}_{2})$ and define the functional
\beginy
\phi _{1}(x_{1})=\int m(x_{1},\underline{x}_{2})f_{2}(\underline{x}_{2})d%
\underline{x}_{2},  \label{int}
\endy
where $f_{2}(\underline{x}_{2})$ is the joint density of $\underline{X}_{2i}$%
. Under the additive structure (\ref{junk}), $\phi _{1}(.)$ is
$m_{1}(.)$ up to a
constant. Replace $m(.)$ in (\ref{int}) with $\hat{\beta}_{0}(x_{1},\underline{x%
}_{2})\equiv\hat{\beta}_{\underline{0}}(\underline{x})$ given by
(\ref{why})
and $\phi _{1}(x_{1})$ can thus be estimated by the sample version of (\ref%
{int}):
\beginn
{\phi}_{n1}(x_{1})=n^{-1}\sum\limits_{i=1}^{n}\hat{\beta}_{0}(x_{1},%
\underline{X}_{2i}).
\endn
As noted by Linton and H\"{a}rdle (1996) and Hengartner and Sperlich
(2005), cautious  choice of the bandwidth is crucial for
${\phi}_{n1}(.)$ to be asymptotically normal. They suggested
different bandwidths  be engaged for the direction
of interest $X_{1}$ and the $d-1$ dimensional nuisance direction $\underline{%
X}_{2}$, say $h_{1}$ and $h$ respectively. Sperlich et al
 (1998) and Linton et al (1999) provide an extensive study of the small sample
properties of marginal integration estimators, including an
evaluation of bandwidth choice.

The following corollary is about the asymptotic properties of
${\phi}_{n1}(.)$.
\begin{Corollary} \label{Cor1} Suppose the support of $\un{X}$ is $[0,1]^{\otimes d}$
with strictly positive density function. Let the conditions in
Proposition \ref{Pro2} hold with $\T_n\equiv\{r(n)/\min(h_1,h)\}^d$
and the $h^d$ replaced by $h_1h^{d-1}$ in all the notations defined
either in (\ref{def}) or (\ref{def1}). If $h_{1}\propto
n^{-1/(2p+3)}$, $h=O(h_{1})$ and (\ref{rbig1}) is modified as
\begin{equation}
\begin{array}{c} nh_1h^{3(d-1)}/\log^3 n \to \infty,\
n^{-1}\{r(n)\}^{\nu_2/2}d_n\log n/M_n^{(2)}\to \infty.
\end{array} \label{rbig0}
\end{equation}
Then we have
\beginn
(nh_{1})^{1/2}\{\phi_{n1}(x_{1})-\phi_{1}(x_{1})\}\overset{L}{\to}%
N(e_{1}W_{p}S_{p}^{-1}B_{1}E\mathbf{m}_{p+1}(x_{1},\underline{X}_{2}),
\tilde\sigma^{2}(x_{1})),\label{davis}
\endn
where `$\stackrel{L}\to$' stands for convergence in distribution,
$$ \tilde \sigma^2(x_1) =\Big\{\int_{[0,1]^{{\otimes d-1}}}
\{fg^2\}^{-1}(x_1,\un{X}_{2})f_2^2(\un{X}_{2})\sigma^2(x_1,\un{X}_{2})d\un{X}_{2}\Big\}e_1S_{p}^{-1}K_2K_2^\t
S_{p}^{-1}e_1^\t ,$$
$\sigma^2(\un{x})=E[\varphi^2(\varepsilon)|\un{X}=\un{x}]$ and $
K_2=\int_{[0,1]^{\otimes d}}K(\un{v})\mu(\un{v})d\un{v} $. In
particular for additive quantile regression, i.e.
$\rho(y;\theta)=(2q-1)(y-\theta)+|y-\theta|$, we have
\beginn
\tilde \sigma^2(x_1) =q(1-q)\Big\{\int_{[0,1]^{{\otimes d-1}}}
f^{-1}(x_1,\un{X}_{2})f^{-2}_{\varepsilon}(0|x_1,\un{X}_{2})
f_2^2(\un{X}_{2})d\un{X}_{2}\Big\}e_1S_{p}^{-1}K_2K_2^\t
S_{p}^{-1}e_1^\t .
\endn
\end{Corollary}

\noindent{\bf Remark 5}. For  conditions in Corollary \ref{Cor1} to
hold, we would need $3d<2p+5$, i.e. the order of local polynomial
approximation increases as the dimension of the covariates
$\underline{X}$ increases. See also the discussion in Hengartner and
Sperlich (2005).

\noindent{\bf Remark 6}. Besides asymptotic normality,  we could
also
 by applying Theorem \ref{T1} develop Bahadur representations for
${\phi } _{n1}(x_{1})$, like those assumed in Linton, Sperlich and
Van Keilegom (2007). Based on (\ref{gexp}), similar results
are also applicable to the generalized additive M-regression model where $G(m(x_{1},\ldots ,x_{d}))$ $%
=c+m_{1}(x_{1})+\ldots +m_{d}(x_{d})$ for some known smooth function
$G(.)$,   in which case the marginal integration estimator is given
by the sample average of $G(\hat{m}(x_{1},\underline{X}_{2i})).$

\section{Concluding Remarks}

Our results can be useful in a variety of contexts including
estimation of quite general nonlinear functionals of M-regression
functions, and we have shown in one specific application how they
can be applied.

\renewcommand{\theLemma}{A\arabic{Lemma}}

 \vspace{.2cm} \section*{Appendix: Regularity
Conditions and Proofs}
 \noindent For any $M>2$, $\lambda_2\in (0,1)$
and $\lambda_1\in (\lambda_2,(1+\lambda_2)/2]$, define
\beginy
&&d_n=(nh^d/\log n)^{-(\lambda_1+{\lambda_2}/{2})}(nh^d\log
n)^{1/2},\ r(n)= (nh^d/\log
n)^{(1-\lambda_2)/2},\label{def}\\
&\n& M_n^{(1)}=M(nh^d/\log n)^{-\lambda_1}, \
M_n^{(2)}=M^{1/4}(nh^d/\log n)^{-\lambda_2}, \
 \T_n=\{r(n)/h\}^d\label{def1}
\endy and $\L_n$ as the smallest
integer such that $\log n(M/2)^{\scriptsize
\L_n+1}>nM_{n}^{(2)}/d_n$. Let  $\|.\|$ denote the Euclidean norm
and $C$ be a generic constant, which may have different values at
each appearance.  Let $\e_i\equiv Y_{i}-m(\un{X}_i)$ and assume that
the following conditions hold.

\begin{description}
\item{(A1)} For each $y\in \R,\ \rho(y;\theta)$ is absolutely
continuous in $\theta,\ i.e.$, there is a function
$\varphi(y;\theta)\equiv\varphi(y-\theta)$ such that for any
$\theta\in \R,\
\rho(y;\theta)=\rho(y;0)+\int^\theta_0\varphi(y;t)dt. $ The
probability density function of $\e_i$ is bounded, $
E\{\varphi(\e_i)|\un{X}_i\}=0 $ almost surely and
$E|\varphi(\e_i)|^{\nu_1}<\infty$ for some $\nu_1>2.$

\item{(A2)} $\varphi(.)$ satisfies the
Lipschitz condition  in $(a_j,a_{j+1}),\ j=0,\cdots,m$, where $a_1
<\cdots < a_m$ are the finite number of jump discontinuity points of
$\varphi(.)$,  $a_0\equiv-\infty$ and $a_{m+1}\equiv+\infty$.

\item{(A3)} $K(.)$ has a compact support, say $[-1,1]^{\otimes d}$ and
$|H_{\un{j}}(\un{u})-H_{\un{j}}(\un{v})|\le C\|u-v\|$ for all $j$
with $0\le |\un{j}|\le 2p+1,$ where
$H_{\un{j}}(u)=\un{u}^{\un{j}}K(\un{u})$.

\item{(A4)} The probability density function of $\un{X},\ f(.)$
is bounded and with bounded first order derivatives. The joint
probability density of $(\un{X}_0,\un{X}_{l})$ satisfies
$f(\un{u},\un{v};l)\le C<\infty$ for all $l\ge 1.$

\item{(A5)} For $\un{r}$ with $|\un{r}|=p+1$, $D^{\un{r}}m(\un{x})$
 is bounded  with bounded first order derivative.

\item{(A6)} The bandwidth $h\to 0$ with
\beginy
nh^d/\log n \to \infty,\ nh^{d+(p+1)/\lambda_2}/\log n<\infty,\ \
n^{-1}\{r(n)\}^{\nu_2/2}d_n\log n/M_n^{(2)}\to \infty,\label{rbig1}
\endy
for some $2<\nu_2\le \nu_1$ and the processes $\{(Y_{i},\un{X}_i)\}$
are strongly mixing with mixing coefficient $\gamma[k]$ satisfying
 \beginy
 &&\sum_{k=1}^\infty k^a\{\gamma[k]\}^{1-2/\nu_2}<\infty \mbox{ for some
} a>(p+d+1)(1-2/\nu_2)/d.\label{centraal}
\endy
Moreover, the bandwidth $h$ and  $\gamma[k]$ should jointly satisfy
the following condition
\beginy
\sum\limits_{n=1}^\infty
n^{3/2}\T_n\Big\{\frac{M_n^{(1)}}{d_n}\Big\}^{1/2}\frac{\gamma[r(n)(2^{\nu_2/2}/M)^{2\L_n/\nu_2}]}
{r(n)(2^{\nu_2/2}/M)^{2\L_n/\nu_2}}\{4M^{2N}\}^{\L_n}<\infty,\
\forall M>0.\label{rbig2}
\endy

\item{(A7)} The conditional density $f_{\un{X}|Y}$ of
$\un{X}$ given $Y$ exists and is bounded. The conditional density
$f_{(\un{X}_1,\un{X}_{l+1})|(Y_1,Y_{l+1})}$ of
$(\un{X}_1,\un{X}_{l+1})$ given $(Y_1,Y_{l+1})$ exists and is
bounded, for all $l\ge 1$.

\end{description}

\noindent{\bf Remark 7}. Assumptions on $\varphi(.)$ in (A1) and
(A2) are satisfied in almost all known robust and likelihood type
regressions. For example, in $qth-$quantile regression, we have
$\varphi(t)=2qI\{t\ge 0\}+(2q-2)I\{t< 0\}$, while for the Huber's
function (\ref{Huber}), its piecewise derivative is given by
\beginn
 \varphi(t)=tI\{|t|<k\}+\mbox{sign}(t)kI\{|t|\ge k\}.
\endn
Note that the condition $E\{\varphi(\e_i)|\un{X}_i\}=0\  a.e.$ is
needed for model specification.  Moreover, if the conditional
density $f(y|\un{x})$ of $Y$ given $\un{X}$ is also continuously
differentiable with respect to $y$, then as proved in Hong (2003)
 there is a constant $C>0$, such that for all small
$t$ and $\un{x}$,
\beginy
E\Big[\Big\{\varphi(Y;t+a)-\varphi(Y;a)\Big\}^2|\un{X}=\un{u}\Big]\le
C|t|\label{A81}
\endy
holds for all  $(a,\un{u})$ in a neighborhood of
$(m(\un{x}),\un{x})$. Define
\beginy
G(t,\un{u})=E\{\varphi(Y;t)|\un{X}=\un{u}\},\quad
G_i(t,\un{u})=(\partial^i/\partial t^i)G(t,\un{u}),\
i=1,2,\label{aa2}
\endy
then it holds that
\beginy g(\un{x})=G_1(m(\un{x}),\un{x})\ge C>0,\ G_2(t,\un{x})
\mbox{ bounded  for all } \un{x}\in {\mathcal D}\mbox{ and }t \mbox{
near }m(\un{x}).\label{aa3}
 \endy
 Assumptions (A3)-(A7) are standard for nonparametric smoothing
in multivariate time series analysis, see Masry (1996). For example,
condition (\ref{centraal}) is needed to bound the covariance of
partial sums of time series as in Lemma \ref{Var1}, while
(\ref{rbig2}) plays a similar role as (4.7b) in Masry (1996). It
guarantees that  the dependence of the time series  is weakly enough
such that the difference caused by the approximation of dependent
random variables by independent ones (through Bradley's strong
approximation theorem) is negligible; see Lemma \ref{mm}. Of course,
(\ref{rbig2}) is more stringent than (4.7b) in Masry (1996), which
is due to the fact that the loss function $\rho(.)$ considered here
is more general than the straightforward square loss.

\noindent{\bf Proof of Proposition \ref{Pro1}}. Write
$\beta^*_{n}(\un{x})=-W_pS_{n,p}^{-1}(\un{x})\sum_{i=1}^{n}Z_{ni}(\un{x})/n$,
where
\beginn
Z_{ni}(\un{x})=H_n^{-1}h^{-d}
K_h(\un{X}_i-\un{x})\varphi(Y_{i},\mu(\un{X}_i-\un{x})^\t
\beta_p(\un{x}))\mu(\un{X}_i-\un{x}).
\endn
We first focus on $EZ_{ni}(\un{x})$. Based on (\ref{aa2}) and
(\ref{aa3}), we have
\beginn
E\{\varphi(Y_{i},\mu(\un{X}_i-\un{x})^\t
\beta_p(\un{x}))|\un{X}_i\}&=&G(\mu(\un{X}_i-\un{x})^\t
\beta_p(\un{x}),\un{X}_i)\\
&=&-g(\un{X}_i)\{m(\un{X}_i)-\mu(\un{X}_i-\un{x})^\t
\beta_p(\un{x})\}\\
&&+G_2(\xi_i(x),\un{X}_i)\{m(\un{X}_i)-\mu(\un{X}_i-\un{x})^\t
\beta_p(\un{x})\}^2/2
\endn
for some $\xi_i(x)$ between $\mu(\un{X}_i-\un{x})^\t
\beta_p(\un{x})$ and $m(\un{X}_i)$. Apparently, if
$\un{X}_i=\un{x}+h\un{v}$, then
\beginn
m(\un{X}_i)-\mu(\un{X}_i-\un{x})^\t
\beta_p(\un{x})=h^{p+1}\sum\limits_{|\un{k}|=p+1}
\frac{D^{\un{r}}m(\un{x})}{\un{k}!}\un{v}^{\un{k}}+h^{p+2}\sum\limits_{|\un{k}|=p+2}
\frac{D^{\un{r}}m(\un{x})}{\un{k}!}\un{v}^{\un{k}}+o(h^{p+2}).
\endn
Therefore,
\beginn
EZ_{ni}(\un{x})&=&h^{p+1}\int K(\un{v})f
g(\un{x}+h\un{v})\mu(\un{v})\sum\limits_{|\un{k}|=p+1}
\frac{D^{\un{r}}m(\un{x})}{\un{k}!}\un{v}^{\un{k}}d\un{v}\\
&&+h^{p+2}\int K(\un{v})f
g(\un{x}+h\un{v})\mu(\un{v})\sum\limits_{|\un{k}|=p+2}
\frac{D^{\un{r}}m(\un{x})}{\un{k}!}\un{v}^{\un{k}}d\un{v}+o(h^{p+2})\\
&\equiv& T_1+T_2.
\endn
Now arrange the $N_{p+1}$ elements of the derivatives
$D^{\un{r}}m(\un{x})/\un{r}!$ for $|\un{r}|=p+1$ as a column vector
${\bf m}_{p+1}(\un{x})$ using the lexicographical order introduced
earlier and define ${\bf m}_{p+2}(\un{x})$ in the similar way. Let
the $N\times N_{p+1}$ matrix $B_{n1}$ and the $N\times N_{p+2}$
matrix $B_{n2}$ be defined as
\beginn
B_{n1}(\un{x})=\left[\begin{array}{c}
                  S_{n,0,p+1}(\un{x})\\
                  S_{n,1,p+1}(\un{x})\\
                  \vdots\\
                  S_{n,p,p+1}(\un{x})
                  \end{array}
      \right],
      \quad
B_{n2}(\un{x})=\left[\begin{array}{c}
                  S_{n,0,p+2}(\un{x})\\
                  S_{n,1,p+2}(\un{x})\\
                  \vdots\\
                  S_{n,p,p+2}(\un{x})
                  \end{array}
      \right],
\endn
where $S_{n,i,p+1}(\un{x})$ and $S_{n,i,p+2}(\un{x})$ is as given by
(\ref{von}). Therefore, $T_1=h^{p+1}B_{n1}(\un{x}){\bf
m}_{p+1}(\un{x}),$ $T_2=h^{p+2}B_{n2}(\un{x}){\bf m}_{p+2}(\un{x}),$
and
\beginn
E\beta^*_{n}(\un{x})=-W_ph^{p+1}S_{n,p}^{-1}(\un{x})B_{n1}(\un{x}){\bf
m}_{p+1}(\un{x})-W_ph^{p+2}S_{n,p}^{-1}(\un{x})B_{n2}(\un{x}){\bf
m}_{p+2}(\un{x})+o(h^{p+2}).
\endn
Let $\un{e}_i,\ i=1,\cdots,d$ be the $d\times 1$ vector having $1$
in the $i$th entry and all other entries $0.$  For $0\le j\le p,\
0\le k\le p+1$, let $N_{j,k}(\un{x})$ be the $N_j\times N_k$ matrix
with its $(l,m)$ element given by
\beginy
\Big[N_{j,k}(\un{x})\Big]_{l,m}=\sum\limits_{i=1}^dD^{\un{e}_i}\{f
g\}(\un{x})\int
K(\un{u})\un{u}^{\tau_{j}(l)+\tau_{k}(m)+\un{e}_i}d\un{u},
\endy
and use these $N_{j,k}(\un{x})$  to construct a $N\times N$ matrix
$N_p(\un{x})$ and a $N\times N_{p+1}$ matrix $\tilde M(\un{x})$
 via
 \beginn
N_p(\un{x})=\left[
                 \begin{array}{cccc}
                  N_{ 0,0}(\un{x})&N_{ 0,1}(\un{x})&\cdots& N_{ 0,p}(\un{x})\\
                  N_{ 1,0}(\un{x})&N_{ 1,1}(\un{x})&\cdots& N_{ 1,p}(\un{x})\\
                  \vdots&\ddots&\vdots\\
                  N_{ p,0}(\un{x})&N_{ p,1}(\un{x})&\cdots& N_{ p,p}(\un{x})\\
                 \end{array}
               \right],\quad
\tilde M(\un{x})=\left[
                 \begin{array}{c}
                  N_{ 0,p+1}(\un{x})\\
                  N_{ 1,p+1}(\un{x})\\
                  \vdots\\
                  N_{ p,p+1}(\un{x})\\
                 \end{array}
               \right].
 \endn
Then $S_{n,p}(\un{x})=\{f g\}(\un{x})S_{p}+hN_{p}(\un{x})+O(h^2),\
B_{n1}(\un{x})=\{f g\}(\un{x})B_{1}+h\tilde M(\un{x})+O(h^2)$ and
$B_{n2}(\un{x})=\{f g\}(\un{x})B_{2}+O(h)$. As
$S_{n,p}^{-1}(\un{x})=\{f g\}^{-1}(\un{x})S_{p}^{-1}-h\{f
g\}^{-2}(\un{x})S_{p}^{-1}N_p(\un{x})S_{p}^{-1}+O(h^2)$, we have
\begin{align*}
-E\beta^*_{n}(\un{x})=&W_ph^{p+1}\Big[\{f
g\}^{-1}(\un{x})S_{p}^{-1}-h\{f
g\}^{-2}(\un{x})S_{p}^{-1}N_p(\un{x})S_{p}^{-1}\Big]\Big[\{f
g\}(\un{x})B_{1}+h\tilde M(\un{x})\Big]{\bf
m}_{p+1}(\un{x})\\
&+W_ph^{p+2}\{f g\}^{-1}(\un{x})S_{p}^{-1}\{f
g\}(\un{x})B_{2}{\bf m}_{p+2}(\un{x})+o(h^{p+2})\\
=&h^{p+1}W_pS_{p}^{-1}B_{1}{\bf m}_{p+1}(\un{x})
+h^{p+2}W_pS_{p}^{-1}\Big[\{f g\}^{-1}(\un{x}){\bf
m}_{p+1}(\un{x})\{\tilde
M(\un{x})-N_p(\un{x})S_{p}^{-1}B_{1}\}\\
& +B_{2}{\bf m}_{p+2}(\un{x})\Big]+o(h^{p+2}).
\end{align*}
We claim that for  elements $E\beta^*_{n\un{r}}(\un{x})$ of
$E\beta^*_{n}(\un{x})$ with $p-|\un{r}|$ even, the $h^{p+1}$ term
will vanish. This means
 for any given $\un{r}$ with $|\un{r}|\le p$ and  $\un{r}_2$ with
 $|\un{r}_2|=p+1$,
\beginy
\sum\limits_{ 0\le |\un{r}|\le
p}\{S_{p}^{-1}\}_{N(\un{r}_1),N(\un{r})}\
\nu_{\un{r}+\un{r}_2}=0.\label{15}
\endy
 To prove this, first note that  for any $\un{r}_1$
with $0\le |\un{r}_1|\le p$ and $\un{r}_2$ with $|\un{r}_2|=p+1$,
\beginy
\sum\limits_{ 0\le |\un{r}|\le
p}\{S_{p}^{-1}\}_{N(\un{r}_1),N(\un{r})}\
\nu_{\un{r}+\un{r}_2}=\int
\un{u}^{\un{r}_2}K_{\un{r}_1,p}(\un{u})d\un{u}, \label{det}
\endy
where
$K_{\un{r},p}(\un{u})=\{|M_{\un{r},p}(\un{u})|/|S_{p}|\}K(\un{u})$
 and $M_{\un{r},p}(\un{u})$ is the same as $S_{p}$, but with the
 $N(\un{r})$ column replaced by $\mu(\un{u})$. Let $c_{ij}$ denote the cofactor of
 $\{S_{p}\}_{i,j}$ and expand the determinant of
 $M_{\un{r},p}(\un{u})$ along the  $N(\un{r})$ column. We see that
 \beginn
\int \un{u}^{\un{r}_2}K_{\un{r},p}(\un{u})d\un{u}=|S_{p}|^{-1}\int
\sum\limits_{ 0\le |\un{r}|\le
p}c_{N(\un{r}),N(\un{r}_1)}\un{u}^{\un{r}_2+\un{r}}K(\un{u})d\un{u}.
 \endn
 (\ref{det}) thus follows, because
$c_{N(\un{r}),N(\un{r}_1)}/|S_{p}|
 =\{S_{p}^{-1}\}_{N(\un{r}_1),N(\un{r})}$ from the symmetry of $S_{p}$
 and a standard result concerning cofactors.
 As a generalization of Lemma 4 in Fan et al (1995) to multivariate case,
 we can further show that for any
 $\un{r}_1$ with $0\le |\un{r}_1|\le p$ and $p-|\un{r}_1|$ even,
 \beginn
\int \un{u}^{\un{r}_2}K_{\un{r},p}(\un{u})d\un{u}=0, \mbox{ for any
}|\un{r}_2|=p+1,
 \endn
 which together with (\ref{det}) yields to
 (\ref{15}).\hspace{\fill}$\Box$

We proceed to prove the main results Theorem \ref{T1}. Define
$\un{X}_{ix}=\un{X}_i-\un{x},\ \mu_{ix}=\mu(\un{X}_{ix}),\
K_{ix}=K_{h}(\un{X}_{ix})$ and
$\varphi_{ni}(\un{x};t)=\varphi(Y_{i};\mu_{ix}^\t\beta_p(\un{x})+t)$.
For $\alpha,\ \beta\in \R^{N}$, define
\beginn
\Phi_{ni}(\un{x};\alpha,\beta)&=&K_{ix}\Big\{\rho(Y_{i};\mu_{ix}^\t(\alpha+\beta+\beta_p(\un{x})))
-\rho(Y_{i};\mu_{ix}^\t(\beta+\beta_p(\un{x})))-\varphi_{i}(\un{x};0)\mu_{ix}^\t\alpha\Big\}\\
&=&K_{ix}\int_{\mu_{ix}^\t\beta}^{\mu_{ix}^\t(\alpha+\beta)}\{\varphi_{ni}(\un{x};t)-\varphi_{ni}(\un{x};0)\}dt,
\endn
and
$R_{ni}(\un{x};\alpha,\beta)=\Phi_{ni}(\un{x};\alpha,\beta)-E\Phi_{ni}(\un{x};\alpha,\beta).$

\begin{Lemma}\label{L2}
Under assumptions $(A1)-(A6)$, we have for all large $M>0,$
\beginy
\sup\limits_{\un{x}\in {\mathcal D}}\sup\limits_{\scriptsize\begin{matrix}\alpha\in B_n^{(1)},\\
\beta\in B_n^{(2)}\end{matrix}}
|\sum\limits_{i=1}^{n}R_{ni}(\un{x};\alpha,\beta)|\le M^{3/2}d_n\
\mbox{almost surely},\label{RR}
\endy
where $\ B_n^{(i)}=\{\beta\in \R^{N}:|H_n\beta|\le M_n^{(i)}\},\
i=1,2.$
\end{Lemma}

\noindent{\bf Proof}. Since ${\mathcal D}$ is compact, it can be
covered by a finite number $\T_n$ of cubes ${\mathcal D}_k={\mathcal
D}_{n,k}$ with side length $l_n=O(\T_n^{-1/d})=O\{h({nh^d}/{\log
n})^{-(1-\lambda_2)/2}\}$ and centers $\un{x}_k=\un{x}_{n,k}$. Write
\begin{align*}
\sup\limits_{\un{x}\in {\mathcal D}}\sup\limits_{\scriptsize\begin{matrix}\alpha\in B_n^{(1)},\\
\beta\in B_n^{(2)}\end{matrix}}
|\sum\limits_{i=1}^{n}R_{ni}(\un{x};\alpha,\beta)|\le&
\max\limits_{\small 1\le k\le \rm\T_n}
\sup\limits_{\scriptsize\begin{matrix}\alpha\in B_n^{(1)},\\
\beta\in B_n^{(2)}\end{matrix}}
\Big|\sum\limits_{i=1}^{n}\Phi_{ni}(\un{x}_k;\alpha,\beta)-E\Phi_{ni}(\un{x}_k;\alpha,\beta)\Big|\\
&+\max\limits_{1\le k\le \T_n}\sup\limits_{\un{x}\in {\mathcal D}_k}
\sup\limits_{\scriptsize\begin{matrix}\alpha\in B_n^{(1)},\\
\beta\in B_n^{(2)}\end{matrix}}
\Big|\sum\limits_{i=1}^{n}\Big\{\Phi_{ni}(\un{x}_k;\alpha,\beta)-\Phi_{ni}(\un{x};\alpha,\beta)\Big\}\Big|\\
&+\max\limits_{1\le k\le \T_n}\sup\limits_{\un{x}\in {\mathcal D}_k}
\sup\limits_{\scriptsize\begin{matrix}\alpha\in B_n^{(1)},\\
\beta\in B_n^{(2)}\end{matrix}}
\Big|\sum\limits_{i=1}^{n}\Big\{E\Phi_{ni}(\un{x}_k;\alpha,\beta)-E\Phi_{ni}(\un{x};\alpha,\beta)\Big\}\Big|\\
\equiv&Q_1+Q_2+Q_3.
\end{align*}
In Lemma \ref{xi}, it is shown that $Q_2\le M^{3/2}d_n/3$ almost
surely and thus $Q_3\le M^{3/2}d_n/3$.

Now all we need to do is to quantify $Q_1$. To this end, we
partition $B_n^{(i)},\ i=1,2,$ into a sequence of disjoint
subrectangles $D_1^{(i)},\cdots,D_{J_1}^{(i)}$ such that
\beginn
|D_{j_1}^{(i)}|=\sup\Big\{|H_{n}(\alpha-\beta)|:\alpha,\beta\in
D_{j_1}^{(i)}\Big\}\le 2M^{-1}M_{n}^{(i)}/\log n, \ \ 1\le j_1\le
J_1.
\endn
Obviously $J_1\le (M\log n)^{N}.$ Choose a point $\alpha_{j_1}\in
D_{j_1}^{(1)}$ and $\beta_{k_1}\in D_{k_1}^{(2)}$. Then
\beginy
\n Q_1&\le &\max\limits_{\scriptsize\begin{matrix}1\le k\le \T_n\\
 1\le j_1,k_1\le J_1\end{matrix}}
\sup\limits_{\scriptsize\begin{matrix}\alpha\in D_{j_1}^{(1)},\\
\beta\in D_{k_1}^{(2)}\end{matrix}}|\sum\limits_{i=1}^{n}
\{R_{ni}(\un{x}_k;\alpha_{j_1},\beta_{k_1})-R_{ni}(\un{x}_k;\alpha,\beta)\}|\\
&&+\max\limits_{\scriptsize\begin{matrix}1\le k\le \T_n\\
 1\le j_1,k_1\le J_1\end{matrix}}
|\sum\limits_{i=1}^{n}R_{ni}(\un{x}_k;\alpha_{j_1},\beta_{k_1})|=H_{n1}+H_{n2}.\label{check}
\endy

We first consider $H_{n1}$. For each $j_1=1,\cdots,J_1$ and $i=1,2$,
partition each rectangle $D_{j_1}^{(i)}$ further into a sequence of
subrectangles $D_{j_1,1}^{(i)},\cdots,D_{j_1,J_2}^{(i)}$. Repeat
this process recursively as follows. Suppose after the $l$th round,
we get a sequence of rectangles $D_{j_1,j_2,\cdots,j_l}^{(i)}$ with
$1\le j_k\le J_k,\ 1\le k\le l$, then in the $(l+1)$th round, each
rectangle $D_{j_1,j_2,\cdots,j_l}^{(i)}$ is partitioned into a
sequence of subrectangles
$\{D_{j_1,j_2,\cdots,j_l,j_{l+1}}^{(i)},1\le j_{l}\le J_{l}\}$ such
that
\beginn
|D_{j_1,j_2,\cdots,j_l,j_{l+1}}^{(i)}|=\sup\Big\{|H_{n}(\alpha-\beta)|:\alpha,\beta\in
D_{j_1,j_2,\cdots,j_l,j_{l+1}}^{(i)}\Big\}\le 2M_{n}^{(i)}/(M^{l}
\log n),\ 1\le j_{l+1}\le J_{l+1},
\endn
where $J_{l+1}\le M^{N}$. End this process after the $(\L_n+1)$th
round, with $\L_n$ given at the beginning of Section 3.  Let
$D_{l}^{(i)},\ i=1,2$, denote the set of all subrectangles of
$D_{0}^{(i)}$ after the $l$th round of partition and a typical
element $D_{j_1,j_2,\cdots,j_l}^{(i)}$ of $D_{l}^{(i)}$ is denoted
as $D_{(j_l)}^{(i)}$. Choose a point $\alpha_{(j_l)}\in
D_{(j_l)}^{(1)}$ and $\beta_{(j_l)}\in D_{(j_l)}^{(2)}$ and define
\beginn
&&\hspace{-.3cm}V_l=\sum\limits_{\tiny\begin{matrix}(j_{l}),\\
(k_{l})\end{matrix}}P\Big\{\Big|\sum\limits_{i=1}^{n}
\{R_{ni}(\un{x}_k;\alpha_{j_l},\beta_{k_l})-R_{ni}(\un{x}_k;\alpha_{j_{l+1}},\beta_{k_{l+1}})\}\Big|
\ge \frac{M^{3/2}d_n}{2^{l}}\Big\},\ 1\le l\le \L_n,\\
&&\hspace{-.3cm}Q_l=\sum\limits_{\tiny\begin{matrix}(j_{l}),\\
(k_{l})\end{matrix}}
P\Big\{\sup\limits_{\tiny\begin{matrix}\alpha\in D_{(j_l)}^{(1)},\\
\beta\in D_{(k_l)}^{(2)}\end{matrix}}\Big|\sum\limits_{i=1}^{n}
\{R_{ni}(\un{x}_k;\alpha_{j_l},\beta_{k_l})-R_{ni}(\un{x}_k;\alpha,\beta)\}\Big|\ge
\frac{M^{3/2}d_n}{2^l}\Big\}, \ 1\le l\le \L_n+1.
\endn
By (A4), it is easy to see that for any $\alpha\in
D_{(j_{\tiny\L_n+1})}^{(1)}\in D_{\tiny\L_n+1}^{(1)}$ and
$\beta\in D_{(k_{\tiny\L_n+1})}^{(2)}\in D_{\tiny\L_n+1}^{(2)}$,
\beginn
|R_{ni}(\un{x}_k;\alpha,\beta_)-R_{ni}(\un{x}_k;\alpha_{j_{\tiny\L_n+1}},\beta_{k_{\tiny\L_n+1}})|\le
\frac{CM_{n}^{(2)}}{ M^{\tiny\L_n+1}\log n},
\endn
which together with the choice of $\L_n$ implies that
$Q_{\tiny\L_n+1}=0.$ As $Q_{l}\le V_l+Q_{l},\ 1\le l\le \L_n,$
\beginy
P(H_{n1}>\frac{M^{3/2}d_n}{2})\le \T_nQ_1\le
\T_n\sum\limits_{l=1}^{\L_n} V_l.\label{Borel}
\endy
To quantify $V_l$, let
 \beginy
W_n=\sum_{i=1}^n Z_{ni},\ Z_{ni}\equiv
R_{ni}(\un{x}_k;\alpha_{j_l},\beta_{k_l})-R_{ni}(\un{x}_k;\alpha_{j_{l+1}},\beta_{j_{l+1}}).\label{exp}
 \endy
Note that by (A2), we have, uniformly in $\un{x},\ \alpha$ and
$\beta$, that
\beginy
|\Phi_{ni}(\un{x};\alpha,\beta)|\le CM_{n}^{(1)}. \label{zzz}
\endy
Therefore, $|Z_{ni}|\le CM_{n}^{(1)}$. With  Lemma \ref{Z}, we can
apply Lemma \ref{mm} to $V_l$ with
\beginn
&&B_1=C_1M_{n}^{(1)},\ B_2=nh^d(M_n^{(1)})^2M_n^{(2)}
\{M^{l}\log n\}^{-2/\nu_2},\\
&& r_n=r_n^l\equiv (2^{\nu_2/2}/M)^{2l/\nu_2}r(n),\ q=n/r_n^l,\ \eta=M^{3/2}d_n/2^{l},\\
&& \lambda_n=(2C_1M_{n}^{(1)}r_n^l)^{-1},\ \Psi(n)=C
q^{3/2}/\eta^{1/2}\gamma[r_n^l]\{r_n^lM_n^{(1)}\}^{1/2}.
\endn
Note that $nM_n^{(1)}/{\eta}\to \infty$, $r_n^l\to \infty$ for all
$1\le l\le \L_n$ from (\ref{rbig1}) and
 \beginn
\lambda\eta=CM^{1/2}\log n M^{2l/\nu_2} /2^{2l},\ \lambda^2B_2=C\log
n^{1-2/\nu_2}M^{2l/\nu_2} /2^{2l}=o(\lambda\eta),
\endn
which  hold uniformly for all $1\le l\le \L_n$. Therefore,
\beginn
V_l\le \Big(\prod\limits_{j=1}^{l+1}J_j^2\Big)4\exp\{-C_1\log n
(M/2^{\nu_2})^{2l/\nu_2}\} +C_2\tau_n^l,
\endn
where, as $J_1\le 2(M\log n )^{N}$ and $J_{l}\le 2M^{N}$ for $2\le
l\le L_n,$ $\tau_n^l$ is given by
\beginn
\tau_n^l= 4^{l}M^{2N(l+1)}(\log
n)^{2N}n^{3/2}\frac{\gamma[r_n^l]\{M_n^{(1)}\}^{1/2}}{r_n^l\{d_n\}^{1/2}}.
\endn
It is tedious but easy to check  that for  $M$ large enough,
\beginy
 \T_n\sum\limits_{l=1}^{\small\L_n}\Big[\Big(\prod\limits_{j=1}^{l+1}J_j^2\Big)4\exp\{-C_1\log n
 (M/2^{\nu_2})^{2l/\nu_2}\}
 \Big]\mbox{ is summable over } n.\label{summa}
\endy
As $\gamma[r_n^l]/r_n^l$ is increasing in $l$, we have
\beginn
\T_n\sum\limits_{l=1}^{\small\L_n}\tau_n^l\le \T_n(\log
n)^{2N}n^{3/2}\frac{\{M_n^{(1)}\}^{1/2}}{\{d_n\}^{1/2}}
\frac{\gamma[r_n^{\tiny\L_n}]}{r_n^{\tiny\L_n}}
\prod\limits_{l=1}^{\small\L_n}4^{l}M^{2N(l+1)},
\endn
which is again summable over $n$ according to (\ref{rbig2}). This
along with (\ref{Borel}) and (\ref{summa})  implies that $H_{n1}\le
{M^{3/2}d_n}/{2}\ $ almost surely, by the Borel-Cantelli lemma.

For $H_{n2}$, first note that
 \beginy
 P(H_{n2}>\eta)
 &\le& \T_nJ_1^2 P(|\sum\limits_{i=1}^{n}R_{ni}(\un{x};\alpha_{j_1},\beta_{k_1})|>\eta).\label{wen}
\endy
We apply Lemma \ref{mm} to quantify
$P(|\sum_{i=1}^{n}R_{ni}(\un{x};\alpha_{j_1},\beta_{k_1}|>\eta)$,
with $r_n=r(n),\   B_1=2C_1M_n^{(1)}$,
$B_2=C_2nh^d(M_n^{(1)})^2M_n^{(2)},\
\lambda_n=\{r(n)M_n^{(1)}\}^{-1}/4C_1$ and $\eta=M^{3/2}d_n$. Then
$nB_1/\eta\to\infty$ and
\beginn
&&\lambda_n\eta/4={(nh^d)^{(1-\lambda_2)/2}(\log n)^{(1+\lambda_2)/2}}/\{16C_1r(n)\}=M^{1/2}\log n/(16C_1),\\
&&\lambda^2_nB_2=M^{1/4}(nh^d)^{1-\lambda_2}(\log n)^{\lambda_2}/\{16C^2_1r^2(n)\}=M^{1/4}\log n/(16C_1^2),\\
&&\Psi(n)\equiv
q_n\{nB_1/\eta\}^{1/2}\gamma[r_n]=\T_nJ_1^2q(n)^{3/2}/\eta^{1/2}\gamma[r(n)]\{r(n)M_n^{(1)}\}^{1/2},
\endn
where $\Psi(n)$ is summable over $n$ by condition (\ref{rbig2}).
Therefore,
\beginy
P(H_{n2}>\eta)\le 2\T_nJ_1^2/n^b+ \Psi(n),\
b=\frac{1}{16C_1}(M^{1/2}-M^{1/4}\frac{C_2}{C_1}).\label{plus}
\endy
By selecting $M$ large enough, we can ensure that (\ref{plus}) is
summable. Thus, for $M$ large enough, $H_{n2}\le M^{3/2}d_n$ almost
surely. By (\ref{check}), we know  for large $M$, $Q_1\le
M^{3/2}d_n$ almost surely.\hspace{\fill}$\Box$

 The quantification of $Q_2$ is very involved, so
we put it as a separate Lemma.
\begin{Lemma}\label{xi}Under the conditions in Lemma \ref{L2}, $Q_2\le M^{3/2}d_n/3$ almost surely.
\end{Lemma}

\noindent{\bf Proof}. Let $\un{X}_{ik}=\un{X}_i-\un{x}_k,\
\mu_{ik}=\mu(\un{X}_{ik}) $ and $ K_{ik}=K_h(\un{X}_{ik}).$ It is
easy to see that  we can write
$\Phi_{ni}(\un{x}_k;\alpha,\beta)-\Phi_{ni}(x;\alpha,\beta)=\xi_{i1}+\xi_{i2}+\xi_{i3}$,
where
\beginn
&&\xi_{i1}=\Big(K_{ik}\mu_{ik}-
K_{ix}\mu_{ix}\Big)^\t\alpha\int_0^1\left\{\varphi_{ni}(\un{x}_k;\mu_{ik}^\t(\beta+\alpha
t))
-\varphi_{ni}(\un{x}_k;0)\right\}dt,\\
&&\xi_{i2}=K_{ix}\mu_{ix}^\t\alpha\int_0^1\left\{\varphi_{ni}(\un{x}_k;\mu_{ik}^\t
(\beta+\alpha t))-\varphi_{ni}(x;\mu_{ix}^\t(\beta+\alpha t))\right\}dt,\\
&&\xi_{i3}=K_{ix}\mu_{ix}^\t\alpha\{\varphi_{ni}(x;0)-\varphi_{ni}(\un{x}_k;0)\}.
\endn
Then $P(Q_2> M^{3/2}d_n/3)\le \T_n(P_{n1}+P_{n2}+P_{n3})$, where
\beginn
P_{nj}\equiv \max\limits_{1\le k\le
\T_n}P\Big(\sup\limits_{\un{x}\in {\mathcal D}_k}
\sup\limits_{\scriptsize\begin{matrix}\alpha\in B_n^{(1)},\\
\beta\in B_n^{(2)}\end{matrix}} |\sum\limits_{i=1}^{n}\xi_{ij}|\ge
{M^{3/2}d_n}/{9}\Big),\ j=1,2,3.
\endn
Based on Borel-Cantelli lemma, $Q_2\le M^{3/2}d_n $ almost surely,
if  $\sum_{n}\T_nP_{nj}<\infty,\ j=1,2,3$.

We first tudy $P_{n1}$. For any fixed $ \alpha\in B_n^{(1)}$ and $
\beta\in B_n^{(2)}$, let $I^{\alpha,\beta}_{ik}=1$, if there exists
some  $t\in [0,1]$, such that there are discontinuity points of
$\varphi(Y_{i};\theta)$ between
$\mu_{ik}^\t(\beta_p(\un{x}_k)+\beta+\alpha t))$ and
$\mu_{ik}^\t\beta_p(\un{x}_k)$; and $I^{\alpha,\beta}_{ik}=0$,
otherwise. Write
$\xi_{i1}=\xi_{i1}I^{\alpha,\beta}_{ik}+\xi_{i1}(1-I^{\alpha,\beta}_{ik})$.
 Note that by (A3), $|(K_{ik}\mu_{ik}-
K_{ix}\mu_{ix})^\t\alpha|\le C_2M_n^{(1)}l_n/h$. Then by (A2) and
the fact that $|\mu_{ik}^\t(\beta+\alpha t)|\le CM_n^{(2)}$, we have
$|\xi_{i1}(1-I^{\alpha,\beta}_{ik})|\le CM_n^{(2)}M_n^{(1)}l_n/h$
uniformly in $i,\alpha$, $\beta$ and $\un{x}\in {\mathcal D}_k$.
Define $U_{ik}=I\{|\un{X}_{ik}|\le 2h\}$, whence
$\xi_{i1}=\xi_{i1}U_{ik}$ since $l_n=o(h)$. Therefore,
\beginy
\n P\Big(\sup\limits_{\scriptsize\begin{matrix}\alpha\in B_n^{(1)},\\
\beta\in B_n^{(2)}\end{matrix}}\sup\limits_{\un{x}\in {\mathcal
D}_k}\Big|\sum\limits_{i=1}^{n}\xi_{i1}(1-I^{\alpha,\beta}_{ik})\Big|>\frac{M^{3/2}d_n}{18}\Big)&\le&
P\Big( \sum\limits_{i=1}^{n}U_{ik}>\frac{M^{1/4}nh^d}{18C}\Big)\\
&\le& P\Big(
|\sum\limits_{i=1}^{n}U_{ik}-EU_{ik}|>\frac{M^{1/4}nh^d}{36C}\Big),
\label{mul}
\endy
where the second inequality follows from the fact that
$\Var(\sum_{i=1}^{n}I\{|\un{X}_{ik}|\le 2h)=O(nh^d)$ implied by
Lemma \ref{Var1}. To quantify (\ref{mul}), we  apply Lemma \ref{mm}
 with $B_1=1, \ \eta=M^{1/4}nh^d/(18C),\ B_2=nh^{d},\ r_n=r(n).$
As  $\lambda_n\eta=CM^{1/4}\log n(nh^d/\log n)^{(1+\lambda_2)/2}$,
$\lambda_n^2B_2=o(\lambda_n\eta)$ and $\T_n\Psi_n$ is summable over
$n$ under condition (\ref{rbig2}), we know that
\beginy
\T_nP\Big(\sup\limits_{\scriptsize\begin{matrix}\alpha\in B_n^{(1)},\\
\beta\in
B_n^{(2)}\end{matrix}}\Big|\sum\limits_{i=1}^{n}\xi_{i1}(1-I^{\alpha,\beta}_{ik})\Big|>M^{3/2}d_n/18\Big)
\mbox{ is summable over }n,\label{zz}
\endy
whence  $\sum_{n}\T_nP_{n1}<\infty,$ is equivalent to
\beginy
\T_nP\Big(\sup\limits_{\scriptsize\begin{matrix}\alpha\in B_n^{(1)},\\
\beta\in
B_n^{(2)}\end{matrix}}\Big|\sum\limits_{i=1}^{n}\xi_{i1}I^{\alpha,\beta}_{ik}\Big|>M^{3/2}d_n/18\Big)
\mbox{ is summable over }n.\label{ggg}
\endy
To prove (\ref{ggg}), first note that $I^{\alpha,\beta}_{ik}\le
I\{\e_i\in S^{\alpha,\beta}_{i;k}\}$, where
 \beginn
S^{\alpha,\beta}_{i;k}&=&\bigcup\limits_{j=1}^m\bigcup\limits_{t\in
[0,1]}[a_j-A(\un{X}_i,\un{x}_k)
+\mu_{ik}^\t(\beta+\alpha t),a_j-A(\un{X}_i,\un{x}_k)]\\
&\subseteq& \bigcup\limits_{j=1}^m
[a_j-CM_n^{(2)},a_j+CM_n^{(2)}]\equiv D_n,\ \mbox{ for some }C>0,\label{shen}\\
A(\un{x}_1,\un{x}_2)&=&(p+1)\sum\limits_{|\un{r}|=p+1}\frac{1}{\un{r}!}(\un{x}_1-\un{x}_2)^{\un{r}}\int_0^1
D^{\un{r}}m(\un{x}_2+w(\un{x}_1-\un{x}_2))(1-w)^pdw,
\endn
where in the derivation of $S^{\alpha,\beta}_{i;k}\subseteq D_n$, we
have used the fact that $|\un{X}_{ik}|\le 2h$ and
$A(\un{X}_i,\un{x}_k)=O(h^{p+1})=O(M_n^{(2)})$ uniformly in $i$. As
$I^{\alpha,\beta}_{ik}\le I\{\e_i\in D_n\}$, we have
$|\xi_{i1}|I^{\alpha,\beta}_{ik}\le |\xi_{i1}| U_{ni}, $ where
$U_{ni}\equiv  I(|\un{X}_{ik}|\le 2h)I\{\e_i\in D_n\}$, which is
independent of the choice of $\alpha$ and $\beta.$
 Therefore,
\begin{align}
\n P\Big(\sup\limits_{\scriptsize\begin{matrix}\alpha\in B_n^{(1)},\\
\beta\in
B_n^{(2)}\end{matrix}}\Big|\sum\limits_{i=1}^{n}\xi_{i1}I^{\alpha,\beta}_{ik}\Big|>M^{3/2}d_n/18\Big)&\le
 P\Big(\sum\limits_{i=1}^{n}
U_{ni}>M^{1/2}nh^dM_n^{(2)}/(18C)\Big)\\
&\le
P\Big(\sum\limits_{i=1}^{n}(U_{ni}-EU_{ni})>\frac{M^{1/2}nh^dM_n^{(2)}}{36C}\Big),\label{ff}
\end{align}
where  the first inequality is because $|\xi_{i1}|\le
CM_n^{(1)}l_n/h$ and the second one is because $EU_{ni}=O(
h^dM_n^{(2)})$ by (A1).
 As $EU_{ni}^2=EU_{ni}$, by Lemma \ref{Var1}, we know that $\Var(\sum_{i=1}^{n}
U_{ni})=Cnh^dM_n^{(2)}$. We can then apply Lemma \ref{mm} to the
last term in   (\ref{ff}) with
\beginn
B_2= Cnh^dM_n^{(2)},\ B_1\equiv 1,\ r_n=r(n),\ \eta\equiv
M^{1/2}nh^dM_n^{(2)}/(36C).
\endn
Apparently, $\lambda_n\eta=C\log n ({nh^d}/{\log
n})^{(1-\lambda_2)/2}$ and $\lambda_n^2B_2=o(\lambda_n\eta).$ As in
this case $\T_n\Psi_n$  is still  summable over $n$ by
(\ref{rbig2}), (\ref{ggg})  thus follows.

For $P_{n2}$, first note that using approach for $P_{n1}$, we can
show that
\beginn
\T_n P\Big(
\sup\limits_{\scriptsize\begin{matrix}\alpha\in B_{n}^{(1)},\\
\beta\in B_{n}^{(2)}\end{matrix}}\sup\limits_{\un{x}\in {\mathcal
D}_k}\Big|\sum_{i=1}^{n}\{\xi_{i2}-\tilde\xi_{i2}\}\Big|\ge
M^{3/2}d_n/18\Big)\mbox{ is summable over }n.
 \endn where
 \beginn
 \tilde\xi_{i2}=K_{ik}\mu_{ik}^\t\alpha\int_0^1\left\{\varphi_{ni}(\un{x}_k;\mu_{ik}^\t(\beta+\alpha t))
 -\varphi_{ni}(x;\mu_{ix}^\t(\beta+\alpha t))\right\}dt.
 \endn
  Therefore, we would have $\sum \T_nP_{n2}<\infty,$ if
 \beginy
 \T_n
P\Big(
\sup\limits_{\scriptsize\begin{matrix}\alpha\in B_{n}^{(1)},\\
\beta\in B_{n}^{(2)}\end{matrix}}\sup\limits_{\un{x}\in {\mathcal
D}_k}\Big|\sum\limits_{i=1}^{n}\tilde\xi_{i2}\Big|\ge
M^{3/2}d_n/18\Big)\mbox{ is summable over }n.\label{han}
 \endy
For any fixed $\alpha\in B_n^{(1)}$, $\beta\in B_n^{(2)}$ and
$\un{x}\in {\mathcal D}_k$, let $I^{\alpha,\beta}_{i;k,x}=1$, if
 there exists some
interval $[t_1,t_2]\subseteq [0,1]$,   such that
\beginy
Y_{i}-\mu_{ik}^\t(\beta_p(\un{x}_k)+\beta+\alpha t)\le a_j\le
Y_{i}-\mu_{ix}^\t(\beta_p(\un{x})+\beta+\alpha t), \ \forall t\in
[t_1,t_2]\label{dis}
\endy
with $a_j\in \{a_1,\cdots,a_m\}$; and $I^{\alpha,\beta}_{i;k,x}=0$,
otherwise. Write
$\tilde\xi_{i2}=\tilde\xi_{i2}I^{\alpha,\beta}_{i;k,x}+\tilde\xi_{i2}(1-I^{\alpha,\beta}_{i;k,x})$.
Note that $K_{ik}\mu_{ik}^\t\alpha=O(M_n^{(1)})$ and
$\varphi_{ni}(\un{x}_k;\mu_{ik}^\t(\beta+\alpha t))
 -\varphi_{ni}(x;\mu_{ix}^\t(\beta+\alpha t))=O(M_n^{(2)}l_n/h)$ if
 $I^{\alpha,\beta}_{i;k,x}=0$. Then again as $\tilde\xi_{i2}=\tilde\xi_{i2}I\{|\un{X}_{ik}|\le 2h\}$,
  we have similar to (\ref{zz})
 that
\beginn
\T_nP\Big(\sup\limits_{\scriptsize\begin{matrix}\alpha\in B_n^{(1)},\\
\beta\in
B_n^{(2)}\end{matrix}}\Big|\sum\limits_{i=1}^{n}\tilde\xi_{i2}(1-I^{\alpha,\beta}_{i;k,x})\Big|>M^{3/2}d_n/18\Big)
\mbox{ is summable over }n.
\endn
Therefore, by (\ref{han}), to show $\sum \T_nP_{n2}<\infty,$ it is
sufficient to show that
\beginy
\T_n P\Big(
\sup\limits_{\scriptsize\begin{matrix}\alpha\in B_{n}^{(1)},\\
\beta\in B_{n}^{(2)}\end{matrix}}\sup\limits_{\un{x}\in {\mathcal
D}_k}\Big|\sum\limits_{i=1}^{n}\tilde\xi_{i2}I^{\alpha,\beta}_{i;k,x}\Big|\ge
M^{3/2}d_n/36\Big)\mbox{ is summable over }n.\label{bund}
\endy
To this end, define $\epsilon_i=\e_i+A(\un{X}_i,\un{x}_k)$. Then
$I^{\alpha,\beta}_{i;k,x}=1$, i.e.  (\ref{dis}) is equivalent to
\beginy
A(\un{X}_i,\un{x}_k)-A(\un{X}_i,\un{x})+\mu_{ix}^\t(\beta+\alpha
t)\le\epsilon_i-a_j\le \mu_{ik}^\t(\beta+\alpha t), \ \forall t\in
[t_1,t_2].\label{dust1}
\endy
 Let $\delta_n\equiv M_n^{(2)}l_n/h$. Then $|A(\un{X}_i,\un{x}_k)-A(\un{X}_i,\un{x})|\le
C\delta_n$, $|(\mu_{ik}-\mu_{ix})^\t\beta| \le C\delta_n$ and
(\ref{dust1}) thus implies that
\beginy
-2C\delta_n+ \mu_{ik}^\t(\beta+\alpha t)\le\epsilon_i-a_j\le
\mu_{ik}^\t(\beta+\alpha t)+2C\delta_n,\
 \ \forall t\in [t_1,t_2].\label{dust2}
\endy
 Without loss of generality, assume $\mu_{ik}^\t\alpha>0.$ Then  from (\ref{dust2})
 we can see that
\beginy
-2C\delta_n+ \mu_{ik}^\t(\beta+\alpha t_2)\le\epsilon_i-a_j\le
\mu_{ik}^\t(\beta+\alpha t_1)+2C\delta_n,\label{dust3}
\endy
which in turn means that if $I^{\alpha,\beta}_{i;k,x}=1$, then
$|\xi_{i2}|\le C(t_2-t_1)|\mu_{ik}^\t\alpha|\le 4C\delta_n$
uniformly in $i,\ \alpha\in B_n^{(1)}$, $\beta\in B_n^{(2)}$ and
$\un{x}\in {\mathcal D}_k$. Therefore, as
$\tilde\xi_{i2}=\tilde\xi_{i2}I\{|\un{X}_{ik}|\le 2h\}$, we have
\beginy
\n && P\Big(
\sup\limits_{\scriptsize\begin{matrix}\alpha\in B_{n}^{(1)}\\
\beta\in B_{n}^{(2)}\end{matrix}}\sup\limits_{\un{x}\in {\mathcal
D}_k}\Big|\sum\limits_{i=1}^{n}\tilde\xi_{i2}I^{\alpha,\beta}_{i;k,x}\Big|\ge
\frac{M^{3/2}d_n}{36}\Big) \\
&& \hspace{.3cm} \le  P\Big(
\sup\limits_{\scriptsize\begin{matrix}\alpha\in B_{n}^{(1)}\\
\beta\in B_{n}^{(2)}\end{matrix}}\sup\limits_{\un{x}\in {\mathcal
D}_k}\sum\limits_{i=1}^{n}I\{|\un{X}_{ik}|\le
2h\}I^{\alpha,\beta}_{i;k,x}\ge
\frac{M^{5/4}nh^dM_n^{(1)}}{36C}\Big). \qquad \label{shan0}
\endy
We will bound $I^{\alpha,\beta}_{i;k,x}$ by a random variable that
is independent of the choice of $\alpha\in B_{n}^{(1)}$ and
$\un{x}\in D_k.$ By the definition of $I^{\alpha,\beta}_{i;k,x}$ and
(\ref{dust3}), the necessary condition for
$I^{\alpha,\beta}_{i;k,x}=1$ is
\beginy
\epsilon_i\in
\bigcup\limits_{j=1}^m[a_j+\mu_{ik}^\t\beta-2M_n^{(1)},a_j+\mu_{ik}^\t\beta+2M_n^{(1)}]\equiv
D_{ni}^{\beta},\label{Je}
\endy
which is indeed independent of the choice of $\alpha$ and $\un{x}\in
{\mathcal D}_k$. Therefore,
\beginy
&&\n P\Big(
\sup\limits_{\scriptsize\begin{matrix}\alpha\in B_{n}^{(1)},\\
\beta\in B_{n}^{(2)}\end{matrix}}\sup\limits_{\un{x}\in {\mathcal
D}_k}\sum\limits_{i=1}^{n}I\{|\un{X}_{ik}|\le
2h\}I^{\alpha,\beta}_{i;k,x}\ge
\frac{M^{5/4}nh^dM_n^{(1)}}{36C}\Big)\\
&\le& P\Big( \sup\limits_{\scriptsize \beta\in
B_{n}^{(2)}}\sum\limits_{i=1}^{n}I\{|\un{X}_{ik}|\le
2h\}I\{\epsilon_i\in D_{ni}^{\beta}\}\ge
\frac{M^{5/4}nh^dM_n^{(1)}}{36C}\Big). \label{shan}
\endy
Now we partition $B_{n}^{(2)}$ into a sequence of subrectangles
$S_1,\cdots,S_{m}$, such that
\beginn
|S_{l}|=\sup\Big\{|H_{n}(\beta-\beta')|:\beta,\beta'\in
S_{l}\Big\}\le M_{n}^{(1)},\ \ 1\le l\le m.
\endn
Obviously, $m\le (M_{n}^{(2)}/M_{n}^{(1)})^N=M^{-3N/4}(nh^d/\log
n)^{(\lambda_1-\lambda_2)N}$.  Choose a point $\beta_l\in S_{l}$ for
each $1\le l\le m,$ and thus
\beginy
\n &&P\Big( \sup\limits_{\scriptsize \beta\in
B_{n}^{(2)}}\sum\limits_{i=1}^{n}I\{|\un{X}_{ik}|\le
2h\}I\{\epsilon_i\in
D_{ni}^{\beta}\}\ge \frac{M^{5/4}nh^dM_n^{(1)}}{36C}\Big)\\
\n &\le& mP\Big( \sum\limits_{i=1}^{n}I\{|\un{X}_{ik}|\le
2h\}I\{\epsilon_i\in D_{ni}^{\beta_l}\}\ge
\frac{M^{5/4}nh^dM_n^{(1)}}{72C}\Big)\\
\n &&+mP\Big( \sup\limits_{\scriptsize \beta'\in
S_{l}}\sum\limits_{i=1}^{n}I\{|\un{X}_{ik}|\le 2h\}|I\{\epsilon_i\in
D_{ni}^{\beta_l}\}-I\{\epsilon_i\in D_{ni}^{\beta'}\}|\ge
\frac{M^{5/4}nh^dM_n^{(1)}}{72C}\Big)\\
&\equiv& m(T_1+T_2). \label{shan1}
\endy
We deal with $T_1$ first. Let
\beginy
U_{ni}^j\equiv I\{|\un{X}_{ik}|\le 2h\}I\{\epsilon_i\in
D_{ni}^{\beta_l}\}.\label{na}
\endy
 Then by
the definition of $D_{ni}^{\beta_j}$ given in (\ref{Je}),
$EU^j_{ni}=O(h^dM_{n}^{(1)}) <M^{5/4}h^dM_n^{(1)}/(144C)$ for large
$M$ and we have
\beginn
T_1\le P\Big( \sum\limits_{i=1}^{n}(U_{ni}^j-EU_{ni}^j)\ge
\frac{M^{5/4}nh^dM_n^{(1)}}{144C}\Big).
\endn
 We can thus apply  Lemma \ref{mm} to the quantity on the
right hand side with $B_1\equiv 1$, $B_2$ given by (\ref{ha0}),
$r_n=r(n)$ and $\eta\propto M^{5/4}nh^dM_n^{(1)}$, and
$\lambda_n=1/(2r_n)$. It follows that
\beginn
\lambda_n\eta=CM^{5/4}\log n(nh^d/\log
n)^{(1+\lambda_2)/2-\lambda_1},\ \lambda_n^2B_2=C\log n(nh^d/\log
n)^{-2(\lambda_1-\lambda_2)/\nu_2}.
\endn
 As $(1+\lambda_2)/2\ge\lambda_1$ and $\lambda_2<\lambda_1$, we have
 $T_1=O( n^{-b})$  for any $b>0.$

For $T_2$, note that as $|\mu_{ik}^\t(\beta-\beta_l)|\le CM_n^{(1)}$
for any $\beta\in S_l$, $1\le l\le m,$ we have
\beginn
|I\{\epsilon_i\in D_{ni}^{\beta_l}\}-I\{\epsilon_i\in
D_{ni}^{\beta}\}|&=&I\{\epsilon_i\in D_{ni}^{\beta_l}\smallsetminus
D_{ni}^{\beta}\}\\
&\le &I\Big\{\epsilon_i\in
\bigcup\limits_{j=1}^m[a_j+\mu_{ik}^\t\beta_l-CM_n^{(1)},a_j+\mu_{ik}^\t\beta_l+CM_n^{(1)}]\Big\}\equiv
U_{ni},
\endn
for some $C>0$, which is independent of the choice of $\beta\in
S_l$. Therefore,
\beginn
T_2\le P\Big( \sum\limits_{i=1}^{n}I\{|\un{X}_{ik}|\le 2h\}U_{ni}\ge
\frac{M^{5/4}nh^dM_n^{(1)}}{72C}\Big),
\endn
which can be dealt with similarly as with $T_1$ and thus $T_2=O(
n^{-b})$ for any $b>0.$ Thus from  (\ref{shan0}), (\ref{shan}) and
(\ref{shan1}), we can claim that (\ref{bund}) is true and thus
$\T_nP_{n2}$ is summable over $n$.

The quantification of $P_{n3}$ is much simpler, as  no $\beta$ is
involved in $\xi_{i3}$. For any given $\un{x}\in {\mathcal D}_k$,
let $I_{i;k,x}=1$, if there is a discontinuity point of
$\varphi(Y_{i};\theta)$ between $\mu_{ik}^\t\beta_p(\un{x}_k)$ and
$\mu_{ix}^\t\beta_p(\un{x})$; and $I_{i;k,x}=0$ otherwise. Write
$\xi_{i3}=\xi_{i3}I_{i;k,x}+\xi_{i3}(1-I_{i;k,x})$. Again by (A2)
and the fact that $|K_{ix}\mu_{ix}^\t\alpha|=O(M_n^{(1)})$ and
$|\mu_{ik}^\t\beta_p(\un{x}_k)-\mu_{ix}^\t\beta_p(\un{x})|=|A(\un{X}_i,\un{x}_k)-A(\un{X}_i,\un{x})|=O(M_n^{(2)}l_n/h)$,
we have similar to (\ref{zz}) that
\beginn
\T_nP\Big(\sup\limits_{\scriptsize\begin{matrix}\alpha\in B_{n}^{(1)}\\
\un{x}\in {\mathcal
D}_k\end{matrix}}\Big|\sum\limits_{i=1}^{n}\xi_{i3}(1-I_{i;k,x})\Big|>M^{3/2}d_n/18\Big)
\mbox{ is summable over }n.
\endn
It's easy to see that $I_{i;k,x}\le I\{\e_i+A(\un{X}_i,\un{x}_k)\in
S_{i;k,x}\}$, where
\beginn
S_{i;k,x}&=&\bigcup\limits_{j=1}^m\bigcup\limits_{t\in
[0,1]}\Big[a_j-|A(\un{X}_i,\un{x}_k)-A(\un{X}_i,\un{x})|,
a_j+|A(\un{X}_i,\un{x}_k)-A(\un{X}_i,\un{x})|\Big]\\
&\subseteq& \bigcup\limits_{j=1}^m
[a_j-CM_n^{(2)}l_n/h,a_j+CM_n^{(2)}l_n/h]\equiv D_n,\mbox{ for some
}C>0.
\endn
Therefore, $|\xi_{i3}|I_{i;k,x}=|\xi_{i3}|I\{|\un{X}_{ik}|\le
2h\}I_{i;k,x}\le U_{ni}$, where
\beginn
U_{ni}\equiv M_n^{(1)}I\{|\un{X}_{ik}|\le
2h\}I\{\e_i+A(\un{X}_i,\un{x}_k) \in D_n\},
\endn
which is independent of the choice of $\alpha\in B_{n}^{(1)}$ and
$\un{x}\in {\mathcal D}_k$. Thus
\beginy
\T_nP\Big(\sup\limits_{\scriptsize\begin{matrix}\alpha\in B_{n}^{(1)}\\
\un{x}\in {\mathcal
D}_k\end{matrix}}\Big|\sum\limits_{i=1}^{n}\xi_{i3}I_{i;k,x}\Big|>M^{3/2}d_n/18\Big)\le
 \T_nP \Big(\sum\limits_{i=1}^{n} [U_{ni}-EU_{ni}]>
M^{3/2}d_n/36\Big ), \label{sarah}
\endy
where we have used the fact that
$EU_{ni}=O(h^dM_n^{(1)}M_n^{(2)}l_n/h)=O(d_n/n)$. We will have $\sum
\T_nP_{n3}<\infty $ if the right hand side in (\ref{sarah}) is
summable over $n$, i.e.
\beginy
 \T_n P \Big(\sum\limits_{i=1}^{n} [U_{ni}-EU_{ni}]> M^{3/2}d_n/36\Big ) \mbox{ is summable over }n.\label{kelly}
\endy
It's easy to check that Lemma \ref{Var1} again holds with
$\psi_{\un{x}}(\un{X}_i,Y_{i})$ standing for $U_{ni}$. Applying
Lemma \ref{mm} to (\ref{kelly}) with $B_1\equiv M_n^{(1)}$,
$B_2\equiv Cnh^d(M_n^{(1)})^2M_n^{(2)}l_n/h$,
 $\eta\equiv M^{3/2}d_n/36$ and $r_n=r(n)$,  we have (note that $nB_1/\eta\to \infty$ indeed)
\beginn
\lambda_n\eta/4=CM^{1/2}\log n,\ \lambda_n^2B_2=Cr_n^{-2/\nu_2}\log
n =o(\lambda_n\eta).
\endn
Thus, $\T_n\Psi_n$ again is summable over $n$ and (\ref{kelly})
indeed holds. \hspace{\fill}$\Box$

\noindent {\bf Proof of Theorem \ref{T1}}. Let
$\lambda_1=\lambda(s)$. Then according to Lemma \ref{L2} and Lemma
\ref{L3}, we know that  with probability $1$, there exists some
$C_1>1$, such that for all large $M>0,$
\beginy
\n &&\sup\limits_{\un{x}\in {\mathcal D}}\sup\limits_{\scriptsize\begin{matrix}\alpha\in B_n^{(1)},\\
\beta\in B_n^{(2)}\end{matrix}}
\Big|\sum\limits_{i=1}^{n}\Phi_{ni}(\un{x};\alpha,\beta)-\frac{nh^d}{2}(H_n\alpha)^\t
S_{np}(\un{x})H_n(\alpha+2\beta)
\Big|\\
&&\le C_1M^{3/2}(d_{n1}+d_n)\le
2C_1M^{3/2}(nh^d)^{1-2\lambda_1}(\log n)^{2\lambda_1}, \mbox{ when }
n \mbox{ is large,}\label{e3}
\endy
where  $d_{n1}=(nh^d)^{1-\lambda_1-2\lambda_2}(\log
n)^{\lambda_1+2\lambda_2}$. Note that from (\ref{matrix}), we can
write
\beginn
\sum\limits_{i=1}^{n}
K_{ni}\varphi(Y_{i};\mu_{ni}^\t\beta_p(\un{x}))\mu_{ni}^\t\alpha=nh^d\beta_n^*(\un{x})^\t
 W_p^{-1}S_{np}(\un{x})H_n\alpha.
\endn
Replace $B_n^{(1)}$ in (\ref{e3}) with $B_{nk}^{(1)}=\Big\{\alpha\in
\R^{N}: k\le M^{-1}(nh^d/\log n)^{\lambda_1}|H_n\alpha|\le k+1\}$
and $M$ with $(k+1)M$.  We have, by the definition of
$\Phi_{ni}(\un{x};\alpha,\beta)$, that
\beginy
\n &\inf\limits_{\un{x}\in {\mathcal D}}&\inf\limits_{\scriptsize\begin{matrix}\alpha\in B_{nk}^{(1)},\\
\beta\in B_n^{(2)}\end{matrix}}\Big\{ \sum\limits_{i=1}^{n}
\rho(Y_{i};\mu_{ni}^\t(\alpha+\beta+\beta_p(\un{x})))K_{ni}
-\sum\limits_{i=1}^{n} \rho(Y_{i};\mu_{ni}^\t(\beta+\beta_p(\un{x})))K_{ni}\\
\n &&+nh^d(W_p^{-1}\beta_n^*(\un{x})-H_n\beta)^\t S_{np}(\un{x})H_n\alpha\Big\}\\
\n &\ge &\inf\limits_{\un{x}\in {\mathcal D}}\inf\limits_{\alpha\in
B_{nk}^{(1)}}\frac{nh^d}{2}(H_n\alpha)^\t S_{np}(\un{x})H_n\alpha-
2CM^{3/2}(nh^d)^{1-2\lambda_1}(\log n)^{2\lambda_1}\\
\n &\ge &\Big\{ C_3(kM)^2/2-2C_1(k+1)^{3/2}M^{3/2}\Big\}(nh^d)^{1-2\lambda_1}(\log n)^{2\lambda_1}\\
&\ge & (8-2^{5/2})C_1C_4^{3/2}(nh^d)^{1-2\lambda_1}(\log
n)^{2\lambda_1}>0 \mbox{ almost surely}, \label{e4}
\endy
where the last term is independent of the choice of $k\ge 1$. The
last inequality is derived as follows. As $S_{p}>0$, suppose its
minimum eigenvalue is $\tau_1>0.$ As $S_{np}(\un{x})\to
g(\un{x})f(\un{x})S_{p}$ uniformly in $\un{x}\in {\mathcal D}$ by
Lemma \ref{Masry}  and $g(\un{x})f(\un{x})$ is bounded away from
zero by (A5) and (\ref{aa3}), there exists some constant $C_3>0$,
such that for all $\un{x}\in {\mathcal D}$, the minimum eigenvalue
of $S_{np}(\un{x})$ is greater than $C_3.$ The last inequality thus
holds if  $M\ge C_4=(16C_1/C_3)^2$. Note that
\beginy
\bigcup\limits_{k=1}^\infty B_{nk}^{(1)}=\Big\{\alpha|\in \R^{N}:
\Big(\frac{nh^d}{\log n}\Big)^{\lambda_1}|H_n\alpha| \ge
M\Big\}:=B_n^{N}.\label{union}
\endy
Therefore, from (\ref{e4}) and (\ref{union}), we have
\beginy
\n &\inf\limits_{\un{x}\in {\mathcal D}}&\inf\limits_{\scriptsize\begin{matrix}\alpha\in B_{n}^{N},\\
\beta\in B_n^{(2)}\end{matrix}}\Big\{ \sum\limits_{i=1}^{n}
\rho(Y_{i};\mu_{ni}^\t(\alpha+\beta+\beta_p(\un{x})))K_{ni}
-\sum\limits_{i=1}^{n} \rho(Y_{i};\mu_{ni}^\t(\beta+\beta_p(\un{x})))K_{ni}\\
 &&+nh^d(W_p^{-1}\beta_n^*(\un{x})-H_n\beta)^\t S_{np}(\un{x})H_n\alpha\Big\}>0\ \mbox{almost surely}.\label{e5}
\endy
Note that by (\ref{fi}), Lemma \ref{L1} and Proposition \ref{Pro1},
we have $|\beta_{n}^*(\un{x})| \le C_3(nh^d/\log n)^{-\lambda_2}$
uniformly in $\un{x}\in {\mathcal D}$ almost surely. Namely,
$\beta_{n}^*(\un{x})\in B_n^{(2)} $
 for all $\un{x}\in {\mathcal D}$, if $M>C_3^4$. This implies that if $M>\max(C_3^4, C_4)$, (\ref{e5}) still holds
 with $\beta$ replaced with $H_n^{-1}W_p^{-1}\beta_n^*(\un{x})$. Therefore,
\beginn
\inf\limits_{\un{x}\in {\mathcal D}}\inf\limits_{\alpha\in
B_{n}^{N}}\Big\{ \sum\limits_{i=1}^{n}
K_{ni}\rho(Y_{i};\mu_{ni}^\t(\alpha
+H_n^{-1}W_p^{-1}\beta_n^*(\un{x})+\beta_p(\un{x})))\hspace{3cm}\\
-\sum\limits_{i=1}^{n}
K_{ni}\rho(Y_{i};\mu_{ni}^\t(H_n^{-1}W_p^{-1}\beta_n^*(\un{x})+\beta_p(\un{x})))\Big\}>0,
\endn
which is equivalent to Theorem \ref{T1}.\hspace{\fill}$\Box$

\noindent{\bf Proof of  (\ref{mathisca})}. Let $\tilde
d_n=(nh^d)^{1-2\lambda_1}(\log n)^{2\lambda_1}$. Through the proof
lines of Theorem \ref{T1},  we can see that (\ref{mathisca}) will
follow if
\beginn
\sup\limits_{\un{x}\in {\mathcal D}}\sup\limits_{\scriptsize\begin{matrix}\alpha\in B_n^{(1)},\\
\beta\in B_n^{(2)}\end{matrix}}
|\sum\limits_{i=1}^{n}R_{ni}(\un{x};\alpha,\beta)|\le M^{3/2}\tilde
d_n\ \mbox{almost surely},\label{RR}
\endn
with $\lambda_1=1,\ \lambda_2=1/2$ and $B_n^{(i)},\ i=1,2$ defined
as in Lemma \ref{L2}.

\noindent To prove this, cover $\mathcal D$ by a finite number
$\tilde \T_n=\{(nh^d/\log n)^{1/2}/h\}^d$ of cubes $\mathcal
D_k=\mathcal D_{nk}$ with side length $\tilde l_n=O\{h(nh^d/\log
n)^{-1/2}\}$ and centers $\un{x}_k=\un{x}_{n,k}$. Write
\begin{align*}
\sup\limits_{\un{x}\in {\mathcal D}}\sup\limits_{\scriptsize\begin{matrix}\alpha\in B_n^{(1)},\\
\beta\in B_n^{(2)}\end{matrix}}
|\sum\limits_{i=1}^{n}R_{ni}(\un{x};\alpha,\beta)|\le&
\max\limits_{\small 1\le k\le \rm\tilde\T_n}
\sup\limits_{\scriptsize\begin{matrix}\alpha\in B_n^{(1)},\\
\beta\in B_n^{(2)}\end{matrix}}
\Big|\sum\limits_{i=1}^{n}\Phi_{ni}(\un{x}_k;\alpha,\beta)-E\Phi_{ni}(\un{x}_k;\alpha,\beta)\Big|\\
&+\max\limits_{1\le k\le \tilde\T_n}\sup\limits_{\un{x}\in {\mathcal
D}_k}
\sup\limits_{\scriptsize\begin{matrix}\alpha\in B_n^{(1)},\\
\beta\in B_n^{(2)}\end{matrix}}
\Big|\sum\limits_{i=1}^{n}\Big\{\Phi_{ni}(\un{x}_k;\alpha,\beta)-\Phi_{ni}(\un{x};\alpha,\beta)\Big\}\Big|\\
&+\max\limits_{1\le k\le \tilde\T_n}\sup\limits_{\un{x}\in {\mathcal
D}_k}
\sup\limits_{\scriptsize\begin{matrix}\alpha\in B_n^{(1)},\\
\beta\in B_n^{(2)}\end{matrix}}
\Big|\sum\limits_{i=1}^{n}\Big\{E\Phi_{ni}(\un{x}_k;\alpha,\beta)-E\Phi_{ni}(\un{x};\alpha,\beta)\Big\}\Big|\\
\equiv&Q_1+Q_2+Q_3.
\end{align*}
We will show that with probability $1$, $Q_k\le M^{3/2}\tilde d_n
/3,\ k=1,2,3$.

\noindent Define $\xi_{ij}$  as in Lemma \ref{L2}. As $P(Q_2>
M^{3/2}\tilde d_n/2)\le \tilde\T_n(P_{n1}+P_{n2}+P_{n3})$, where
\beginn
P_{nj}\equiv \max\limits_{1\le k\le
\tilde\T_n}P\Big(\sup\limits_{\un{x}\in {\mathcal D}_k}
\sup\limits_{\scriptsize\begin{matrix}\alpha\in B_n^{(1)},\\
\beta\in B_n^{(2)}\end{matrix}} |\sum\limits_{i=1}^{n}\xi_{ij}|\ge
{M^{3/2}\tilde d_n}/{9}\Big),\ j=1,2,3.
\endn
Then based on Borel-Cantelli lemma, $Q_2\le M^{3/2}\tilde d_n/2 $
almost surely if $\sum_{n}\tilde\T_nP_{nj}<\infty,$ for $\ j=1,2,3.$
We only prove that for $P_{n1}$ to illustrate. Recall that
\beginn
\xi_{i1}=\Big(K_{ik}\mu_{ik}-
K_{ix}\mu_{ix}\Big)^\t\alpha\int_0^1\left\{\varphi_{ni}(\un{x}_k;\mu_{ik}^\t(\beta+\alpha
t)) -\varphi_{ni}(\un{x}_k;0)\right\}dt.
\endn
Because $|(K_{ik}\mu_{ik}- K_{ix}\mu_{ix})^\t\alpha|\le
C_2M_n^{(1)}\tilde l_n/h$, $|\mu_{ik}^\t(\beta+\alpha t)|\le
CM_n^{(2)}$ and $\varphi(.)$ is Lipschitz continuous, we have
$|\xi_{i1}|\le CM_n^{(2)}M_n^{(1)}\tilde l_n/h$. Define
$U_{ik}=I\{|\un{X}_{ik}|\le 2h\}$. As $\tilde l_n=o(h)$, we can see
that $\xi_{i1}=\xi_{i1}U_{ik}$ and similar to (\ref{mul}), we have
\beginn
\n P\Big(\sup\limits_{\scriptsize\begin{matrix}\alpha\in B_n^{(1)},\\
\beta\in B_n^{(2)}\end{matrix}}\sup\limits_{\un{x}\in {\mathcal
D}_k}\Big|\sum\limits_{i=1}^{n}\xi_{i1}\Big|>\frac{M^{3/2}\tilde
d_n}{9}\Big)&\le&
P\Big( \sum\limits_{i=1}^{n}U_{ik}>\frac{M^{1/4}nh^d}{9C}\Big)\\
&\le& P\Big(
|\sum\limits_{i=1}^{n}U_{ik}-EU_{ik}|>\frac{M^{1/4}nh^d}{18C}\Big),
\endn
and $\sum_{n}\tilde\T_nP_{nj}<\infty$ thus follows from similar
arguments as those lying between (\ref{mul}) and (\ref{zz}).

\noindent The proof of $Q_1\le M^{3/2}\tilde d_n/2 $ almost surely
is much easier than  in Lemma \ref{L2}, if $\varphi(.)$ is Lipschitz
continuous. Instead of the iterative partition approach adopted
there, we once for all partition $B_n^{(i)},\ i=1,2,$ into a
sequence of disjoint subrectangles $D_1^{(i)},\cdots,D_{J_1}^{(i)}$
such that
\beginn
|D_{j_1}^{(i)}|=\sup\Big\{|H_{n}(\alpha-\beta)|:\alpha,\beta\in
D_{j_1}^{(i)}\Big\}\le M_{n}^{(i)}(\log n/n)^{1/2}, \ \ 1\le j_1\le
J_1.
\endn
Obviously $J_1\le (n/\log n)^{N/2}.$ Choose a point $\alpha_{j_1}\in
D_{j_1}^{(1)}$ and $\beta_{k_1}\in D_{k_1}^{(2)}$. Then
\beginy
\n Q_1&\le &\max\limits_{\scriptsize\begin{matrix}1\le k\le \tilde\T_n\\
 1\le j_1,k_1\le J_1\end{matrix}}
\sup\limits_{\scriptsize\begin{matrix}\alpha\in D_{j_1}^{(1)},\\
\beta\in D_{k_1}^{(2)}\end{matrix}}|\sum\limits_{i=1}^{n}
\{R_{ni}(\un{x}_k;\alpha_{j_1},\beta_{k_1})-R_{ni}(\un{x}_k;\alpha,\beta)\}|\\
&&+\max\limits_{\scriptsize\begin{matrix}1\le k\le \T_n\\
 1\le j_1,k_1\le J_1\end{matrix}}
|\sum\limits_{i=1}^{n}R_{ni}(\un{x}_k;\alpha_{j_1},\beta_{k_1})|=H_{n1}+H_{n2}.\label{check}
\endy
By  Lipschitz continuity of $\varphi(.)$, we have for any $\alpha\in
D_{j_1}^{(1)}$ and $ \beta\in D_{k_1}^{(2)}$,
\beginn
|\Phi_{ni}(\un{x}_k;\alpha_{j_1},\beta_{k_1})-\Phi_{ni}(\un{x}_k;\alpha,\beta)|^2=O(\{M_{n}^{(2)}\}^3\log
n/n)<M^{3/2}\tilde d_n/(4n).
\endn
Therefore, it remains to show that $P(H_{n2}>M^{3/2}\tilde d_n/4)$
is summable over $n$.

\noindent First note that by  Cauchy inequality,
$|R_{ni}(\un{x};\alpha,\beta)|^2=O(\{M_n^{(1)}M_n^{(2)}\}^2)$ and
$E|R_{ni}(\un{x};\alpha,\beta)|^2=O(h^d\{M_n^{(1)}M_n^{(2)}\}^2)$
uniformly in $\un{X}_i,\ \un{x},\ \alpha\in M_n^{(1)}$ and $\beta\in
M_n^{(2)}$. Next, for any $\eta>0,$
\beginn
P(H_{n2}>\eta)
 &\le& \tilde\T_nJ_1^2 P(|\sum\limits_{i=1}^{n}R_{ni}(\un{x};\alpha_{j_1},\beta_{k_1})|>\eta).
\endn
We apply Lemma \ref{mm} with $r_n=(nh^d/\log n)^{1/2},\
B_1=2C_1M_n^{(1)}M_n^{(2)}$, $B_2=C_2nh^d(M_n^{(1)}M_n^{(2)})^2,\
\lambda_n=(4C_1r_n\{M_n^{(2)}\}^{2})^{-1}$ and $\eta=M^{3/2}\tilde
d_n/4$. It is easy to see that $nB_1/\eta\to\infty$ and
\beginn
&&\lambda_n\eta/4=M\log n/(16C_1),\ \lambda^2_nB_2=o(\lambda_n\eta)\\
&&\Psi(n)\equiv q_n\{nB_1/\eta\}^{1/2}\gamma[r_n]=n^{3/2}(\log n
)^{-1/2}\gamma[r(n)]/r(n).
\endn
As  $\tilde \T_nJ_1^2\Psi(n)$ is summable over $n$ by condition
(\ref{rbig2}), so is $P(H_{n2}>M^{3/2}\tilde d_n/4)$.
\hspace{\fill}$\Box$

\noindent{\bf Proof of Corollary \ref{Pro2}}. As $1+\lambda_2\ge
2\lambda_1$,
 it's sufficient to prove that with probability $1$,
\begin{equation}
\beta^*_{n}(\un{x})-E\beta^*_{n}(\un{x})-\frac{1}{nh^d}
W_pS_{np}^{-1}(\un{x})H_n^{-1}\sum\limits_{i=1}^{n}
K_h(\un{X}_i-\un{x})\varphi(\varepsilon_{i})\mu(\un{X}_i-\un{x})=
O\Big\{\Big(\frac{\log
n}{nh^d}\Big)^{(1+\lambda_2)/2}\Big\},\label{fi}
\end{equation}
uniformly in $\un{x}\in {\mathcal D}$.
 As $\varphi(\varepsilon_{i})\equiv\varphi(Y_{i},m(X_i))$ and
$E\varphi(\varepsilon_{i})=0$,  the term on the left hand side of
(\ref{fi}) stands for
\beginn
W_pS_{n,p}^{-1}(\un{x})\frac{1}{nh^d}\sum_{i=1}^{n}\{Z_{ni}(\un{x})-EZ_{ni}(\un{x})\},
\endn
where
\beginn
Z_{ni}(\un{x})=H_n^{-1}
K_h(\un{X}_i-\un{x})\mu(\un{X}_i-\un{x})\Big\{\varphi(Y_{i},\mu(\un{X}_i-\un{x})^\t
\beta_p(\un{x}))-\varphi(\varepsilon_{i})\Big\}.
\endn
Next, like what we did in Lemma \ref{L2}, we cover ${\mathcal D}$
with number $\T_n$
 cubes ${\mathcal D}_k={\mathcal D}_{n,k}$ with side length
$l_n=O(\T_n^{-1/d})$ and centers $\un{x}_k=\un{x}_{n,k}$. Write
\begin{align*}
\sup\limits_{\un{x}\in {\mathcal D}}
|\sum\limits_{i=1}^{n}Z_{ni}(\un{x})-EZ_{ni}(\un{x})|\le&
\max\limits_{\small 1\le k\le \rm\T_n}
\Big|\sum\limits_{i=1}^{n}Z_{ni}(\un{x}_k)-EZ_{ni}(\un{x}_k)\Big|\\
&+\max\limits_{1\le k\le \T_n}\sup\limits_{\un{x}\in {\mathcal D}_k}
\Big|\sum\limits_{i=1}^{n}Z_{ni}(\un{x})-Z_{ni}(\un{x}_k)\Big|\\
&+\max\limits_{1\le k\le \T_n}\sup\limits_{\un{x}\in {\mathcal D}_k}
\Big|\sum\limits_{i=1}^{n}EZ_{ni}(\un{x})-EZ_{ni}(\un{x}_k)\Big|\\
\equiv&Q_1+Q_2+Q_3.
\end{align*}
As $Z_{ni}(\un{x})-Z_{ni}(\un{x}_k)=H_n^{-1}
K_h(\un{X}_i-\un{x})\mu(\un{X}_i-\un{x})\{\varphi_{ni}(\un{x};0)-\varphi_{ni}(\un{x}_k;0)\}$,
through approaches similar to that for $\xi_{i3}$ in the proof of
Lemma \ref{xi}, we can show that
\beginn
Q_2=O\Big\{\Big(\frac{nh^d}{\log n}\Big)^{(1-\lambda_2)/2}\log
n\Big\}\mbox{ almost surely}
\endn
 and so is $Q_3.$ To bound $Q_1,$ first note that
$EZ_{ni}^2(\un{x}_k)=O(h^{p+1+d})$ uniformly in $i$ and $k.$ As
$|Z_{ni}(\un{x})|\le C$ for some constant $C$ by (A2), we can see
that from Lemma \ref{Var1}
\beginn
\sum\limits_{i=1}^{n}EZ_{ni}^2(\un{x}_k)+\sum\limits_{i<
j}|\Cov(Z_{ni}(\un{x}_k),Z_{nj}(\un{x}_k))|\le C_2nh^{p+1+d}.
\endn
Finally by Lemma \ref{mm} with $B_1=C_1$ , $B_2\equiv Cnh^{p+1+d}$,
 $\eta=A_3 (nh^d/\log n)^{(1-\lambda_2)/2}\log
n$ and $r_n=r(n)$, we have (note that $nB_1/\eta\to \infty$ indeed)
\beginn
\lambda_n\eta= A_3/(2C_1)\log n,\ \lambda_n^2B_2=C_2/(4C_1^2)\log n.
\endn
Therefore,
\beginn
P\Big(\max\limits_{\small 1\le k\le \rm\T_n}
\Big|\sum\limits_{i=1}^{n}Z_{ni}(\un{x}_k)-EZ_{ni}(\un{x}_k)\Big|\ge
A_3(nh^d/\log n)^{(1-\lambda_2)/2}\log n\Big)\le \T_n/n^a+C
\T_n\Psi_n,
\endn
where $\ a=A_3/(8C_1)-C_2/(4C_1^2)$. By selecting $A_3$ large
enough, we can ensure that $\T_n/n^a$ is summable over $n$. As
$\T_n\Psi_n$ is summable over $n$ from (\ref{rbig2}), we can
conclude that
\beginn
Q_1=O\Big\{\Big(\frac{nh^d}{\log n}\Big)^{(1-\lambda_2)/2}\log
n\Big\}\mbox{ almost surely.}
\endn
This together with Lemma \ref{Masry} completes the
proof.\hspace{\fill}$\Box$

\noindent{\bf Proof of Corollary \ref{Cor1}}. Through the proof
lines for Theorem \ref{T1} and Corollary \ref{Pro2},
 it's not difficult to see that Corollary \ref{Pro2} still holds
 under the conditions imposed here. Under
 the additive structure
(\ref{junk}), we thus have
\begin{align}
\n\phi_{n1}(x_1)=&\phi_1(x_1)+\frac{1}{n}\sum\limits_{i=1}^nm_2(\un{X}_{2i})
-h^{p+1}e_{1}W_pS_{p}^{-1}B_{1}\frac{1}{n}\sum\limits_{i=1}^n{\bf
m}_{p+1}(x_1,\un{X}_{2i})\\
\n&+\frac{1}{n^2h_1h^{d-1}}e_1\sum\limits_{j=1}^n\varphi(\varepsilon_{j})
\sum\limits_{i=1}^nS_{np}^{-1}(x_1,\un{X}_{2i})
K(X_{1,xj}/h_1,\un{X}_{2,ij}/h)\mu(X_{1,xj}/h_1,\un{X}_{2,ij}/h)\\
&+o_p(\{\max(h_1,h)\}^{p+1})+O_p\{({nh_1h^{d-1}}/{\log
n})^{-3/4}\},\label{sunny}
\end{align}
 where $X_{1,xj}=X_{1j}-x$,
$\un{X}_{2,ij}=\un{X}_{2i}-\un{X}_{2j}$
 and $e_1$ is as in Proposition
\ref{Pro1}. Note that by (\ref{rbig0}),
$(nh_1)^{1/2}({nh_1h^{d-1}}/{\log n})^{-3/4}\to 0,$ the $O_p(.)$
term can thus be safely ignored.\\
By  central limit theorem for strongly mixing processes (Bosq, 1998,
Theorem 1.7), we have
\beginn
\frac{1}{n}\sum\limits_{i=1}^nm_2(\un{X}_{2i})=O_p(n^{-1/2}),\quad
\frac{1}{n}\sum\limits_{i=1}^n{\bf m}_{p+1}(x_1,\un{X}_{2i})=E{\bf
m}_{p+1}(x_1,\un{X}_{2})+O_p(n^{-1/2}).
\endn
As the expectations of all other terms in (\ref{sunny}) are  $0$,
the leading term in the asymptotic bias of
$\tilde\phi_1(x_1)-\phi_1(x_1)$ is thus given by
\beginn
-\{\max(h_1,h)\}^{p+1}e_{1}W_pS_{p}^{-1}B_{1}E{\bf
m}_{p+1}(x_1,\un{X}_{2}).
\endn
Again through standard arguments in Masry (1996), we can see that
\beginn
&&\frac{1}{nh^{d-1}}\sum\limits_{i=1}^nS_{np}^{-1}(x_1,\un{X}_{2i})
K_h(X_{1,xj},\un{X}_{2,ij})\mu(X_{1,xj}/h_1,\un{X}_{2,ij}/h)\\
&=&S_{np}^{-1}(x_1,\un{X}_{2j})f_2(\un{X}_{2j})\int_{[0,1]^{\otimes
d-1}}
\{K\mu\}(X_{1,xj}/h_1,\un{v})d\un{v}\Big\{1+O\Big(\Big\{\frac{\log
n}{nh^{d-1}}\Big\}^{1/2}\Big)\Big\}
\endn
uniformly in $1\le i\le n.$ Therefore, the leading term in the
asymptotic variance of $\phi_{n1}(x_1)-\phi_1(x_1)$ is the variance
of the following term
\beginn
&&(nh_1)^{-1}e_1\sum\limits_{j=1}^n\varphi(\varepsilon_{j})
S_{np}^{-1}(x_1,\un{X}_{2j})f_2(\un{X}_{2j})\int_{[0,1]^{\otimes
d-1}} \{K\mu\}(X_{1,xj}/h_1,\un{v})d\un{v},
\endn
which is asymptotically
\beginy
(nh_1)^{-1}\Big\{\int_{[0,1]^{\otimes d-1}}
\{fg^2\}^{-1}(x_1,\un{X}_{2})f_2^2(\un{X}_{2})\sigma^2(x_1,\un{X}_{2})d\un{X}_{2}\Big\}e_1S_{p}^{-1}K_2K_2^\t
S_{p}^{-1}e_1^\t.\label{qq}
\endy
If $\rho(y;\theta)=(2q-1)(y-\theta)+|y-\theta|$ and
$\varphi(\theta)=2qI\{\theta>0\}+(2q-2)I\{\theta<0\}$, we have
$g(\un{x})=2f_{\varepsilon}(0|\un{x})$ and
\beginn
&&\sigma^2(\un{x})=E[\varphi^2(\varepsilon)|\un{X}=\un{x}]=
4q^2(1-F_{\varepsilon}(0))+4(1-q)^2F_{\varepsilon}(0)=4q(1-q),
\endn
which when substituted into (\ref{qq}), yields the asymptotic
variance  for the quantile regression estimator,
\beginn
\tilde \sigma^2(x_1) =q(1-q)\Big\{\int_{[0,1]^{\otimes d-1}}
f^{-1}(x_1,\un{X}_{2})f^{-2}_{\varepsilon}(0|x_1,\un{X}_{2})f_2^2(\un{X}_{2})d\un{X}_{2}\Big\}e_1S_{p}^{-1}K_2K_2^\t
S_{p}^{-1}e_1^\t. \hspace{2cm}\Box
\endn
 The next Lemma  is due to Davydov (Hall and
Heyde (1980), Corollary A.2).
\begin{Lemma}\label{L0}
 Suppose that $X$ and $Y$ are random variables which are
$\mathcal G-$ and $\mathcal H-$ measurable, respectively, and that
$E|X|^p<\infty$, $E|Y|^q<\infty,$ where $p,\ q>1,\ p^{-1}+q^{-1}<1$.
Then
\beginn
|EXY-EXEY|\le 8\|X\|_p\|Y\|_q\Big\{\sup\limits_{A\in \mathcal G,B\in
\mathcal H}|P(AB)-P(A)P(B)|\Big\}^{1-p^{-1}-q^{-1}}.
\endn
\end{Lemma}
The next lemma is a generalization of some results in the proof of
Theorem 2 in Masry (1996).
\begin{Lemma}\label{mm}
Suppose $\{Z_i\}_{i=1}^\infty$ is a zero-mean strictly stationary
processes with strongly mixing coefficient $\gamma[k]$, and that
$|Z_i|\le B_1$, $\sum_{i=1}^nEZ_{i}^2+\sum_{i<
j}|\Cov(Z_{i},Z_{j})|\le  B_2$. Then for any $\eta>0$ and integer
series $r_n\to \infty$, if $nB_1/{\eta}\to \infty$ and
$q_n\equiv[n/r_n]\to \infty,$ we have
\beginn
 P(|\sum\limits_{i=1}^{n} Z_i|\ge \eta)\le
 4\exp\{-\frac{\lambda_n\eta}{4}+\lambda_n^2B_2\}+C\Psi(n),
 \endn
where $\Psi(n)=q_n\{nB_1/\eta\}^{1/2}\gamma[r_n],\ \lambda_n=1/\{2r_nB_1\}$.
\end{Lemma}

\noindent{\bf Proof}. We partition the set $\{1,\cdots,n\}$ into $
2q\equiv 2q_n$ consecutive blocks of size $r\equiv r_n$ with
$n=2qr+v$ and $0\le v<r$. Write
\beginn
V_n(j)= \sum\limits_{i=(j-1)r+1}^{jr}Z_i,\ j=1,\cdots,2q
\endn
and
\beginn
W_n'= \sum\limits_{j=1}^{q}V_n(2j-1),\ W_n''=
\sum\limits_{j=1}^{q}V_n(2j),\ W_n'''=
\sum\limits_{i=2qr+1}^{n}Z_{i}.\label{friend}
\endn
Then $W_n\equiv \sum_{i=1}^n Z_i=W_n'+W_n''+W_n'''$. The
contribution of $W_n'''$ is negligible as it consists of at most
$r$ terms compared of $qr$ terms in $W_n'$ or $W_n''$. Then by the
stationarity of the processes, for any $\eta>0,$
\beginy
P(W_n>\eta)\le P(W_n'>\eta/2)+P(W_n''>\eta/2)=2 P(W_n'>\eta/2).\label{use1}
\endy
To bound $P(W_n'>\eta/2)$, using recursively  Bradley's Lemma, we
can approximate the random variables
$V_n(1),V_n(3),\cdots,V_n(2q-1)$ by independent random variables
$V^*_n(1),V^*_n(3),$ $\cdots,V^*_n(2q-1)$, which satisfy that for
$1\le j\le q$, $V^*_n(2j-1)$ has the same distribution as
$V_n(2j-1)$ and
\beginy
P\Big(|V^*_n(2j-1)-V_n(2j-1)|>u\Big)\le
18({\|V_n(2j-1)\|_\infty}/{u})^{1/2}\sup|P(AB)-P(A)P(B)|,\label{garden}
\endy
where $u$ is any positive value such that $0<u\le
\|V_n(2j-1)\|_\infty<\infty $ and  the supremum is taken over
 all sets of $A$ and $B$  in the $\sigma-$algebras of
events generated by $\{V_n(1),V_n(3),\cdots,V_n(2j-3)\}$ and
$V_n(2j-1)$ respectively. By the definition of $V_n(j)$, we can see
that $\sup|P(AB)-P(A)P(B)|=\gamma[r_n]$.
 Write
\beginy
\n P(W_n'>\frac{\eta}{2})&&\hspace{-.4cm}\le P\Big(\Big|\sum\limits_{j=1}^{q}V^*_n(2j-1)\Big|>\frac{\eta}{4}\Big)
+P\Big(\Big|\sum\limits_{j=1}^{q}V_n(2j-1)-V^*_n(2j-1)\Big|>\frac{\eta}{4}\Big)\\
&&\hspace{-.4cm}\equiv I_1+I_2.\label{sun}
\endy
We bound $I_1$ as follows.  Let $\lambda=1/\{2B_1r\}.$ Since
$|Z_i|\le B_1$, $\lambda|V_n(j)|\le 1/2,$ then using the fact that
$e^x\le 1+x+x^2/2$ holds for $|x|\le 1/2$, we have
\beginy
E\Big\{e^{\pm\lambda V^*_n(2j-1)}\Big\}\le 1+\lambda^2E\{V_n(j)\}^2\le e^{\lambda^2E\{V^*_n(2j-1)\}^2}.\label{al}
\endy
By Markov inequality, (\ref{al}) and the independence of the $\{V^*_n(2j-1)\}_{j=1}^q$, we have
\beginy
\n I_1&\le& e^{-\lambda\eta/4}\Big[E\exp\Big(\lambda\sum\limits_{j=1}^qV^*_n(2j-1)\Big)
+E\exp\Big(-\lambda\sum\limits_{j=1}^qV^*_n(2j-1)\Big)\Big]\\
\n &\le& 2\exp\Big(-\lambda\eta/4+\lambda^2\sum\limits_{j=1}^qE\{V^*_n(2j-1)\}^2\Big)\\
&\le& 2\exp\Big\{-\lambda\eta/4+C_2\lambda^2B_2\Big\}. \label{Al}
\endy
We now bound the term $I_2$ in (\ref{sun}). Notice that
\beginn
I_2\le\sum\limits_{j=1}^q P\Big(\Big|V_n(2j-1)-V^*_n(2j-1)\Big|>\frac{\eta}{4q}\Big).
\endn
If $\|V_n(2j-1)\|_\infty\ge{\eta}/(4q)$,  substitute ${\eta}/(4q)$ for $u$ in (\ref{garden}),
\beginy
I_2\le 18q\{\|V_n(2j-1)\|/{\eta}/(4q)\}^{1/2}\gamma[r_n]\le Cq^{3/2}/\eta^{1/2}\gamma[r_n](r_nB_1)^{1/2},\label{silence}
\endy
If $\|V_n(2j-1)\|_\infty<{\eta}/(4q)$, let $u\equiv \|V_n(2j-1)\|_\infty$ in (\ref{garden}) and we have
\beginn
I_2\le Cq\gamma[r_n],
\endn
which is of smaller order than (\ref{silence}), if $nB_1/{\eta}\to
\infty.$ Thus by (\ref{use1}), (\ref{sun}), (\ref{Al}) and
(\ref{silence}),
\beginn
 P(W_n>\eta)\le 4\exp\{-\lambda_n\eta/4+C_2B_2\lambda^2_n\}
+C\Psi_n,
\endn
where the constant $C$ is independent of $n$.\hspace{\fill}$\Box$

\begin{Lemma}\label{Var1} For any $\un{x}\in R^d$, let $\psi_{\un{x}}(\un{X}_i,Y_{i})=I(|\un{X}_{ix}|\le h)
\psi_x(\un{X}_{ix},Y_{i})$, a  measurable function of $(\un{X}_i,Y_{i})$ with $|\psi_{\un{x}}(\un{X}_i,Y_{i})|\le B$
 and $V=E\psi_{\un{x}}^2(\un{X}_i,Y_{i})$. Suppose the mixing coefficient $\gamma[k]$ satisfies
 (\ref{centraal}). Then
\beginn
\Cov(\sum\limits_{i=1}^{n}|\psi_{\un{x}}(\un{X}_i,Y_{i})|)=
nV\Big[1+o\Big\{\Big(B^2h^{p+d+1}/V\Big)^{1-2/\nu_2}\Big\}\Big].
\endn
\end{Lemma}

\noindent{\bf Proof}. Denote $\psi_{\un{x}}(\un{X}_i,Y_{i})$ by
$\psi_{ix}$. First note that
\beginn
V=E\psi_{ix}^2=h^d\int_{\scriptsize
|\un{u}|\le 1}
E(\psi_{ix}^2|\un{X}_i=\un{x}+h\un{u} )f(\un{x}+h\un{u})d\un{u},
\endn
\beginn
\sum\limits_{i< j}|\Cov(\psi_{ix},\psi_{jx})|&=&\sum\limits_{l=1}^{n-d}(n-l-d+1)|\Cov(\psi_{0x},\psi_{lx})|
\le n\sum\limits_{l=1}^{n-d}|\Cov(\psi_{0x},\psi_{lx})|\\
&=&n\sum\limits_{l=1}^{d-1}+n\sum\limits_{l=d}^{\pi_n}+n\sum\limits_{l=\pi_n+1}^{n-d}\equiv nJ_{21}+nJ_{22}+nJ_{23},
\endn
where $\pi_n=h^{(p+d+1)(2/\nu_2-1)/a}$. For $J_{21}$, there might be
an overlap between the components of $\un{X}_0$ and $\un{X}_l$, for
example, when  $\un{X}_i=(X_{i-d},\cdots,X_{i-1})$, where $\{X_i\}$
is a univariate time series. Without loss of generality,
 let $\un{u}',\un{u}''$ and $\un{u}'''$ of dimensions $l,d-l$ and $l$ respectively, be the  $d+l$ distinct
 random variables in $(\un{X}_{0x}/h,\un{X}_{lx}/h)$. Write $\un{u}_1=(\un{u}'^\t,\un{u}''^\t)^\t$ and $
\un{u}_2=(\un{u}''^\t,\un{u}'''^\t)^\t$. Then  by Cauchy inequality, we have
\beginy
\Big|E\Big(\psi_{0x},\psi_{lx}|
{\scriptsize\begin{matrix}
\un{X}_0=\un{x}+h\un{u}_1\\
\un{X}_l=\un{x}+h\un{u}_2
\end{matrix}}
\Big)\Big|\le \Big\{E(\psi_{0x}^2|\un{X}_0=\un{x}+h\un{u}_1)E(\psi_{jx}^2|\un{X}_j=\un{x}+h\un{u}_2)\Big\}^{1/2}
=V/h^d\label{f}
\endy
and through a transformation of variables, we have
\beginn
|\Cov(\psi_{0x},\psi_{lx})|\le h^{l}V\hspace{-.2cm}\int_{\scriptsize
\begin{matrix}
|\un{u}_1|\le 1\\
|\un{u}_2|\le 1
\end{matrix}}
|f(\un{x}+h\un{u}_1,\un{x}+h\un{u}_2;l)-
f(\un{x}+h\un{u}_1)f(\un{x}+h\un{u}_2;l+d-1)|d\un{u}'d\un{u}''d\un{u}''',
\endn
where by (A4) and (A5), the integral is bounded. Therefore,
\beginn
nJ_{21}\le CnV\sum\limits_{l=1}^{d-1}h^{l}=o(nV).
\endn
For $J_{22}$,  there is no overlap between the components of $\un{X}_0$ and $\un{X}_l$. Let $\un{X}_{0x}=h\un{u}$
and $\un{X}_{lx}=h\un{v}$ and we have
\beginn
|\Cov(\psi_{0x},\psi_{lx})|&\le &h^{2d}\int_{\scriptsize
\begin{matrix}
|\un{u}|\le 1\\
|\un{v}|\le 1
\end{matrix}}  E\Big(\psi_{0x},\psi_{lx}|
{\scriptsize\begin{matrix}
\un{X}_0=\un{x}+h\un{u}\\
\un{X}_l=\un{x}+h\un{v}
\end{matrix}}
\Big)d\un{u}d\un{v}\\
&& \qquad \qquad \times  [f(\un{x}+h\un{u},\un{x}+h\un{v};l+d-1)- f(\un{x}+h\un{u})f(\un{x}+h\un{v})]\\
&=&Ch^{d}V,
\endn
where the last equality  follows from (A4), (A5) and (\ref{f}).
Therefore, as $\pi_nh^{d}\to 0$,
\beginn
nJ_{22}=O\{n\pi_nh^{d}V\}=o(nV).
\endn
For $J_{23}$, using Davydov's lemma (Lemma \ref{L0}) we have
\beginy
|\Cov(\psi_{0x},\psi_{lx})|\le
8\{\gamma[l-d+1]\}^{1-2/\nu_2}\{E|\psi_{ix}|^{\nu_2}\}^{2/\nu_2},\
\mbox{ as }\nu_2>2.\label{alba}
\endy
As $|\psi_{ix}|\le B$, $E|\Phi_{ni}|^{\nu_2}\le B^{\nu_2-2}V$,
\beginn
J_{23}\le CB^{(\nu-2)2/\nu_2}
V^{2/\nu_2}/\pi_n^{a}\sum\limits_{l=\pi_n+1}^{\infty}l^a\{\gamma[l-d+1]\}^{1-2/\nu_2},
\endn
where the summation term is $o(1)$ as $\pi_n\to \infty$. Thus
$J_{23}=o\Big\{V\Big(B^2h^{p+d+1}/V\Big)^{1-2/\nu_2}\Big\}$, which
 completes the proof. \hspace{\fill}$\Box$

\begin{Lemma}\label{Z}
Suppose (A2)- (A6) hold. Then for $U_{ni}^l, l=1,\cdots,m$ defined
in (\ref{na}) and $Z_{ni}, l=1,\cdots,\L_n$  defined in
(\ref{exp}), we have
\beginy
&&\hspace{-.5cm}\sum\limits_{i=1}^{n}E (U_{ni}^l)^2+\sum\limits_{i<
j}|\Cov(U^l_{ni},U^{l}_{nj})|\le Cnh^dM_n^{(1)}\{M_n^{(2)}/M_n^{(1)}\}^{1-2/\nu_2},\label{ha0}\\
&&\hspace{-.5cm} \sum\limits_{i=1}^{n}EZ^2_{ni}+\sum\limits_{i<
j}|\Cov(Z_{ni},Z_{nj})|= nh^d(M_n^{(1)})^2M_n^{(2)}\{M^{l}\log
n\}^{-2/\nu_2},\label{ha2}
\endy uniformly in $\un{x}_k$, $1\le k\le \T_n$.
\end{Lemma}

\noindent{\bf Proof}. We only prove (\ref{ha2}), which is more
involved than (\ref{ha0}). To simplify the notations, denote
$\alpha_{j_l},\beta_{k_l},\alpha_{j_{l}}$ and $\beta_{j_{l}}$ by
$\alpha_1,\beta_1,\alpha_2$ and $\beta_2$, respectively. Clearly,
\begin{equation*}
\int_{\un{u}^\t H_n\beta_2}^{\un{u}^\t
H_n(\alpha_2+\beta_2)}\{\varphi_{ni}(\un{x}_k;t)-\varphi_{ni}(\un{x}_k;0)\}dt
=\int_{\un{u}^\t H_n\beta_1}^{\un{u}^\t
H_n(\alpha_2+\beta_1)}\{\varphi_{ni}(\un{x}_k;t+\un{u}^\t
H_n(\beta_2-\beta_1)) -\varphi_{ni}(\un{x}_k;0)\}dt,
\end{equation*}
and
\beginn
 Z_{ni}&=&\int_{\un{u}^\t H_n\beta_1}^{\un{u}^\t H_n(\alpha_1+\beta_1)}\{\varphi_{ni}(\un{x}_k;t)
 -\varphi_{ni}(\un{x}_k;0)\}dt
 -\int_{\un{u}^\t H_n\beta_2}^{u^\t H_n(\alpha_2+\beta_2)}\{\varphi_{ni}(\un{x}_k;t)-\varphi_{ni}(\un{x}_k;0)\}dt\\
 &=&\int_{\un{u}^\t H_n\beta_1}^{\un{u}^\t H_n(\alpha_1+\beta_1)}\{\varphi_{ni}(\un{x}_k;t)-\varphi_{ni}(\un{x}_k;t+
 \un{u}^\t H_n(\beta_2-\beta_1))\}dt\\
 &&-\int_{\un{u}^\t H_n(\alpha_1+\beta_1)}^{\un{u}^\t H_n(\alpha_2+\beta_1)}\{\varphi_{ni}(\un{x}_k;t+\un{u}^\t
 H_n(\beta_2-\beta_1))-\varphi_{ni}(\un{x}_k;0)\}dt\equiv \Delta_1+\Delta_2.
 \endn
Therefore, $E\{Z_{ni}\}^2=h^d\int K^2(\un{u})f(\un{x}_k+h\un{u})
E\{(\Delta_1+\Delta_2)^2|X_i=\un{x}_k+h\un{u}\}d\un{u}.$
The conclusion is thus obvious observing that by Cauchy inequality and (\ref{A81}),
\beginn
E(\Delta_1^2|X_i=\un{x}_k+h\un{u})&\le &|\un{u}^\t H_n\alpha_1
\un{u}^\t H_n(\beta_2-\beta_1)\un{u}^\t H_n\alpha_1|
\le 2(M_{n}^{(1)})^2M_{n}^{(2)}/(M^{l}\log n),\\
 E(\Delta_2^2|X_i=\un{x}_k+h\un{u})&\le &\{\un{u}^\t H_n(\alpha_2-\alpha_1)\}^2(|\un{u}^\t H_n\alpha_2|
 +|\un{u}^\t H_n\alpha_1|+2|\un{u}^\t H_n\beta_2|)\\
&\le & 4 (M_{n}^{(1)})^2M_{n}^{(2)}/(M^{l}\log n)^2,
 \endn
where we used the facts that $|\alpha_1-\alpha_2|\le
2M_{n}^{(1)}/(M^{l}\log n)$ and $|\beta_1-\beta_2|\le
 2 M_{n}^{(2)}/(M^{l}\log n)$. Therefore, $E\{Z_{ni}\}^2=Ch^d(M_{n}^{(1)})^2M_{n}^{(2)}/(M^{l}\log n)$.
 As $|Z_{ni}|\le CM_{n}^{(1)}$ and $h^{p+1}/M_{n}^{(2)}<\infty$, the rest
 of the proof can be completed following the proof of Lemma \ref{Var1}.
\hspace{\fill}$\Box$

\begin{Lemma}\label{Var}Suppose (A2)- (A6) hold.
\beginy
&&\sum\limits_{i=1}^{n}E\Phi^2_{ni}+\sum\limits_{i<
j}|\Cov(\Phi_{ni},\Phi_{nj})|\le Cnh^d(M_n^{(1)})^2M_n^{(2)},
\label{ha1}
\endy
\end{Lemma}
uniformly in $\un{x}\in {\mathcal D},\alpha\in B_{n}^{(1)}$ and
$\beta\in B_{n}^{(2)}$.

\noindent{\bf Proof}. By Cauchy inequality and (\ref{A81}), we have
\beginy
\n && E\Phi^2_{ni}\\
\n &\hspace{-.3cm}=&\hspace{-.3cm}h^d\int K^2(\un{u})
E\Big[\Big\{\int_{\mu(\un{u})^\t H_n\beta}^{\mu(\un{u})^\t
H_n(\alpha+\beta)}\Big(\varphi_{ni}(\un{x};t)-\varphi_{ni}(\un{x};0)\Big)dt\Big\}^2|\un{X}_i
=\un{x}+h\un{u}\Big]f(\un{x}+h\un{u})d\un{u}\\
\n &\hspace{-.3cm}\le&\hspace{-.3cm}h^d\int
f(\un{x}+h\un{u})K^2(\un{u})\mu(\un{u})^\t H_n\alpha\int_{\un{u}^\t
H_n\beta}^{\mu(\un{u})^\t
H_n(\alpha+\beta)}E\Big[\Big(\varphi_{ni}(\un{x};t)-\varphi_{ni}
(\un{x};0)\Big)^2|\un{X}_i=\un{x}+h\un{u}\Big]dtd\un{u}\\
&\hspace{-.3cm}\le&\hspace{-.3cm}h^d\int K^2(\un{u})\mu(\un{u})^\t
H_n\alpha\int_{\mu(\un{u})^\t H_n\beta}^{\mu(\un{u})^\t
 H_n(\alpha+\beta)}C|t|dtf(\un{x}+h\un{u})d\un{u}=O\Big\{ h^d(M_n^{(1)})^2M_n^{(2)}\Big\},\label{sophie}
\endy
uniformly in $\un{x}\in {\mathcal D}$,  $\alpha\in B_{n}^{(1)}$ and
$\beta\in B_{n}^{(2)}$.  (\ref{ha1}) thus follows from
(\ref{sophie}) and Lemma \ref{Var1}. \hspace{\fill}$\Box$

\begin{Lemma}\label{Masry}Let $(A3)-(A6)$
hold. Then
\beginn
\sup\limits_{\un{x}\in {\mathcal
D}}|S_{np}(\un{x})-g(\un{x})f(\un{x})S_p|=O(h+(nh^d/\log n)^{-1/2} )
\ \mbox{almost surely}.
\endn
\end{Lemma}
{\bf Proof}. The result is almost the same as  Theorem 2 in Masry
(1996). Especailly if (\ref{rbig2}) holds, then the requirement
(3.8a)
 there on the mixing coefficient $\gamma[k]$ is met.\hspace{\fill}$\Box$

\begin{Lemma}\label{L3}
Denote $d_{n1}=(nh^d)^{1-\lambda_1-2\lambda_2}(\log
n)^{\lambda_1+2\lambda_2}$ and let $\lambda_1$ and $B_n^{(i)},\
i=1,2,$ be as in Lemma \ref{L2}. Suppose that $(A1)-(A5)$ and
(\ref{rbig1}) hold. Then there is a constant $C>0$ such that for
each $M>0$ and all large $n,$
\beginn
\sup\limits_{\un{x}\in {\mathcal
D}}\sup\limits_{\scriptsize\begin{matrix}\alpha\in
B_n^{(1)},\\\beta\in B_n^{(2)}\end{matrix}}
|\sum\limits_{i=1}^{n}E\Phi_{ni}(\un{x};\alpha,\beta)-\frac{nh^d}{2}(H_n\alpha)^\t
S_{np}(\un{x})H_n(\alpha+2\beta)|\le
 CM^{3/2}d_{n1}.
\endn
\end{Lemma}
\noindent{\bf Proof}. Recall that
$G(t,\un{u})=E(\varphi(Y;t)|\un{X}=\un{u})$,
\beginy
E\Phi_{ni}(\un{x};\alpha,\beta)&=&h^d\int K(\un{u})f(\un{x}+h\un{u})d\un{u}\times
\int_{\mu(\un{u})^\t H_n\beta}^{\mu(\un{u})^\t H_n(\alpha+\beta)}\label{double}\\
\n &&\Big\{G(t+\mu(\un{u})^\t
H_n\beta_p(\un{x}),\un{x}+h\un{u})-G(\mu(\un{u})^\t
H_n\beta_p(\un{x}),\un{x}+h\un{u})\Big\}dt.
\endy
By $(A3)$ and $(A5)$, we have
\beginn
&&G(t+\mu(\un{u})^\t H_n\beta_p(\un{x}),\un{x}+h\un{u})-G(\mu(\un{u})^\t H_n\beta_p(\un{x}),\un{x}+h\un{u})\\
&& \hspace{3cm}=tG_1(\mu(\un{u})^\t H_n\beta_p(\un{x}),\un{x}+h\un{u})
+\frac{t^2}{2}G_2(\xi_n(t,\un{u};\un{x}),\un{x}+h\un{u}),\\
&&G_1(\mu(\un{u})^\t
H_n\beta_p(\un{x}),\un{x}+h\un{u})=g(\un{x}+h\un{u})+O(h^{p+1}),
\endn
where $\xi_n(t,\un{u};\un{x})$ falls between $\mu(\un{u})^\t
H_n\beta_p(\un{x})$ and $t+\mu(\un{u})^\t H_n\beta_p(\un{x})$, and
the term $O(h^{p+1})$ is uniform  in $\un{x}\in {\mathcal D}$.
Therefore, the inner integral in (\ref{double}) is given by
\beginn
\frac{1}{2}g(\un{x}+h\un{u})(H_n\alpha)^\t\mu(\un{u})\mu(\un{u})^\t
H_n(\alpha+2\beta)+O\Big\{M^{3/2} \Big(\frac{\log
n}{nh^d}\Big)^{\lambda_1+2\lambda_2}\Big\}
\endn
 uniformly in $\un{x}\in {\mathcal D}$, where we have used the fact
 that $nh^{d+(p+1)/\lambda_2}/\log n<\infty$.
 By the definition of $S_{np}(\un{x})$, the proof is thus completed.\hspace{\fill}$\Box$

 \begin{Lemma}\label{L1}Under conditions in Theorem  \ref{T1}, we
 have
 \beginn
\sup\limits_{\un{x}\in {\mathcal D}}\Big|\frac{1}{nh^d}
W_pS_{np}^{-1}(\un{x})H_n^{-1}\sum\limits_{i=1}^{n}
K_h(\un{X}_i-\un{x})\varphi(\varepsilon_{i})\mu(\un{X}_i-\un{x})\Big|=O\Big\{\Big(\frac{\log
n}{nh^d}\Big)^{1/2}\Big\} \ \mbox{almost surely}.
\endn
\end{Lemma}

\noindent{\bf Proof}. Note that, under conditions Theorem  \ref{T1},
the assumptions imposed
  by Masry (1996) in Theorem 5 are validated. Specifically, (4.5) there follows
  from (\ref{rbig1}) and
  (4.7b) there can be derived  from (\ref{rbig2}).
  Therefore, following the proof lines there, we can show that
 \beginn
 \sup\limits_{\un{x}\in {\mathcal D}}\Big|\frac{1}{nh^d}
H_n^{-1}\sum\limits_{i=1}^{n}
K_h(\un{X}_i-\un{x})\varphi(\varepsilon_{i})\mu(\un{X}_i-\un{x})\Big|=O\Big\{\Big(\frac{\log
n}{nh^d}\Big)^{1/2}\Big\},
 \endn
which together with Lemma \ref{Masry} yields the desired
results.\hspace{\fill}$\Box$

%
%
%\renewcommand{\theequation}{A\arabic{equation}}
%\setcounter{equation}{0}
%\begin{center}
% {{\large A}PPENDIX} \\
% {\it Assumptions and Proofs}
% \end{center}
%
%

\

\begin{center}
 {{\large R}EFERENCES}
\end{center}

\begin{description}

\ditem Andrews, D.W.K. (1994).  Asymptotics for semiparametric
econometric models via stochastic equicontinuity. {\it
Econometrica} {\bf 62}, 43-72.

\ditem Bahadur, R.R. (1966). A note on quantiles in large samples.
{\it Ann. Math. Statist}. {\bf 37}, 577-80.

\ditem Bosq, D. (1998). {\it Nonparametric Statistics for Stochastic
Processes}. NewYork: Springer-Verlag.

\ditem Chen, X., Linton, O. B.  and I. Van Keilegom (2003).
Estimation of Semiparametric Models when the Criterion is not
Smooth. {\it Econometrica} {\bf 71}, 1591-608.

\ditem Fan, J., Heckman, N.E. and Wand, M.P. (1995). Local
polynomial kernel regression for generalized linear models and
quasi-likelihood functions. {\it J. Amer. Statist. Assoc.} {\bf 90},
141-50.

\ditem Fan, J. and Gijbels, I. (1996). {\it Local polynomial regression. }%
London: Chapman and Hall.

\ditem Hengartner, N. W. and Sperlich, S. (2005).     Rate optimal
estimation with the integration method in the presence of many
covariates. {\it J.  Multivariate Anal.} {\bf  95}, 246 - 72.

\ditem  Hall, P. and Heyde, C.C. (1980). {\it Martingale Limit
Theory and its Applications}. NewYork: Academic Press.

\ditem Hong, S. (2003).  Bahadur representation and its application
for Local Polynomial Estimates in Nonparametric M-Regression. {\it
J. Nonparametric Statist.} {\bf 15}, 237-51.

\ditem Horowitz, J. L. and Lee, S. (2005). Nonparametric estimation
of an additive quantile regression model. {\it J. Amer. Statist.
Assoc.} {\bf 100},  1238-49.

\ditem Huber, P. J. (1973) Robust regression. {\it Ann. Statist.}
{\bf 1},  799-821.

\ditem Kiefer, J. (1967). On Bahadur's representation of sample
quantiles.{ \it Ann. Math. Statist}. {\bf 38}, 1323-42.

\ditem Linton, O. B. (2001). Estimating additive nonparametric
models by partial $L_q$ Norm: The Curse of Fractionality. {\it
Econom. Theory} {\bf 17}, 1037-50.

\ditem Linton, O. B., Hardle, W and Sperlich, S (1999). A Simulation
comparison between the Backfitting and Integration methods of
estimating Separable Nonparametric Models. {\it TEST} {\bf 8},
419-58.

\ditem Linton, O. B. and Nielsen, J. P. (1995). A kernel method of
estimating structured nonparametric regression based on marginal
integration. {\it Biometrika} {\bf 82}, 93-100.

\ditem Linton, O. B., Sperlich, S.  and I. Van Keilegom (2007).
Estimation of a Semiparametric Transformation Model by Minimum
Distance. {\it  Ann.  Statist.} To appear

\ditem Linton, O. B. and H\"ardle, W. (1996). Estimation of additive
regression models with known links. {\it Biometrika} {\bf 83},
529-40.

 \ditem Masry, E. (1996). Multivariate local polynomial
regression for time series: uniform strong consistency and rates.
{\it J. Time Ser. Anal.} {\bf 17}, 571-99.

\ditem Peng, L. and Yao, Q. (2003). Least absolute deviation
estimation for ARCH and GARCH models. {\it Biometrika} \textbf{90},
967-75.

\ditem Rosenblatt, M. A central limit theorem and strong mixing
conditions. {\it Prof Nat. Acad. Sci.} {\bf 4}, 43-7.

\ditem Sperlich, S., O. Linton, and W. H\"{a}rdle (1998). A
Simulation comparison between the Backfitting and Integration
methods of estimating separable nonparametric models. {\it Test}
{\bf 8}, 419-58.

\ditem Stone, C. J. (1982). Optimal global rates of convergence for
nonparametric regression. {\it Ann. Statist.} {\bf 10}, 1040-53.

\ditem Stone, C. J. (1986). The dimensionality reduction principle
for generalized additive models. {\it Ann. Statist.} {\bf 14},
592-606.

\ditem Wu, W. B. (2005). On the Bahadur representation of sample
quantiles for dependent sequences. {\it Ann. Statist.} {\bf 33},
1934-63.

\end{description}

\end{document}